\documentclass[11pt]{article}

\usepackage[margin=1in]{geometry}
\usepackage{multirow}
\usepackage[T1]{fontenc}
\usepackage[english]{babel}
\usepackage{lmodern}
\usepackage{microtype}
\usepackage{hyphenat} 
\usepackage{subcaption}
\usepackage{enumitem}
\usepackage{authblk} 
\usepackage{tabularx}
\newcolumntype{C}{>{\centering\arraybackslash}X}

\usepackage[authoryear,round]{natbib}

\usepackage[colorlinks=true,citecolor=blue,linkcolor=blue,urlcolor=blue]{hyperref}

\usepackage[dvipsnames]{xcolor}
\usepackage{graphicx}
\usepackage{pgfplots}
\pgfplotsset{compat=1.18}
\usepackage{tikz}
\usetikzlibrary{arrows,automata,positioning,math}

\usepackage{amsthm,amssymb}
\usepackage{mathtools} 
\usepackage{bbm}       
\usepackage{booktabs}

\theoremstyle{plain}
\newtheorem{theorem}{Theorem}[section]
\newtheorem{lemma}[theorem]{Lemma}

\newtheorem{proposition}[theorem]{Proposition}

\theoremstyle{definition}

\theoremstyle{remark}

\numberwithin{equation}{section}

\usepackage{thmtools,thm-restate}

\newcommand{\bP}{\mathbb{P}}
\newcommand{\bE}{\mathbb{E}}
\newcommand{\N}{\mathbb{N}}
\newcommand{\R}{\mathbb{R}}

\newcommand{\1}{\mathbbm{1}}
\newcommand{\cP}{\mathcal{P}}
\usepackage{float}
\usepackage{ifthen}

\definecolor{forestgreen}{rgb}{0.13,0.55,0.13}

\title{Mean-Variance Finite-Menu Pricing}
\author[1]{Tim S.\,G.\,van Eck\thanks{Email: \texttt{t.s.g.vaneck@tilburguniversity.edu}}}
\author[1]{Pieter Kleer}
\author[1]{Johan S.\,H.\,van Leeuwaarden}
\affil[1]{Department of Econometrics and Operations Research, Tilburg University}
\date{}

\begin{document}
\maketitle

\begin{abstract}
How much pricing flexibility is needed when the distribution of buyers' willingness to pay is unknown but its mean and variance are available? The two extremes, a robust posted price and unrestricted price randomization, are well understood, yet the optimal mechanism with an arbitrary finite number of prices has not been characterized. We solve this mean--variance problem for every menu size, recovering optimal prices, randomization probabilities, and a worst-case valuation distribution from a system of two equations in two unknowns. The solution uncovers a reversal of the posted pricing prescription: increasing variance lowers guaranteed revenue but eventually \emph{raises the expected offered price} as soon as the seller can randomize over two or more prices. A low price protects against the downside, while higher prices preserve access to high-valuation buyers. We also study pricing under mean absolute deviation (MAD) information and establish an equivalence with mean--support ambiguity that reduces the problem to one with a known solution.
Unlike under variance, every optimal price decreases as MAD increases, even though MAD permits much heavier tails. Finally, we examine how dispersion affects the revenue gains from offering additional prices.
\end{abstract}

\section{Introduction}
How should a seller set prices when the distribution of buyers' willingness to pay is unknown? A particularly natural starting point is to assume that the mean and variance of valuations can be estimated, even if the demand curve itself cannot. These two statistics are familiar, economically interpretable, and far less demanding than a fully specified valuation model. They also form a classical basis for distributionally robust pricing: the seller selects a policy that performs well against every distribution consistent with the available moments. If the valuation distribution were known, posting a single price would be optimal \citep{myerson1981optimal,riley1983optimal}. However, under distributional ambiguity, the seller may benefit from randomizing over several prices. Against a single posted price, the worst-case distribution places a large share of the probability just below that price, so that many buyers narrowly fail to purchase. A randomized mechanism offers no such single price to target, allowing the seller to hedge against the worst case \citep{carrasco2018optimal}.

The seller's options thus range from a single posted price to a randomization over a continuum of prices, the most general way of selling to a single buyer. Under this type of ambiguity, both extremes are well understood. The robust posted price is known \citep{azar2012optimal,chen2022distribution}, and so is the optimal mechanism that randomizes over a continuum of prices, which generally secures a substantially higher revenue guarantee \citep{carrasco2018optimal,chen2021screening}. Neither is entirely satisfactory, however. Fully randomized mechanisms might be difficult to implement in practice, both from an operational and market-design perspective. A single posted price, on the other hand, may capture only a small fraction of the revenue that full randomization secures. This raises a natural question. \emph{What happens in between, when the seller can randomize over only a small, fixed number of prices?}

\citet{wang2025power} initiates the study of this middle ground by developing a general framework for robust pricing with finite menus under various forms of distributional information, including support, mean, and quantile restrictions. Under a monotonicity condition on the functions defining the ambiguity set, the seller’s problem for given prices reduces to a finite linear program. Wang shows that relatively few prices can capture much of the benefit of unrestricted randomization and, under mean--support information, characterizes optimal menus of every size. Mean--variance information does not satisfy this monotonicity condition and requires a different approach.

As our first contribution, we develop an exact solution to the mean--variance finite-menu pricing problem, jointly determining the optimal prices, their probabilities, and a worst-case distribution. Despite the nonconvexity of the problem, the solution is determined by a system of two equations in two unknowns, irrespective of the menu size.

\begin{figure}[h!]
  \centering
  \begin{subfigure}[b]{0.49\linewidth}
    \centering
    \caption*{Lognormal ground-truth distribution}
    \includegraphics[width=\linewidth]{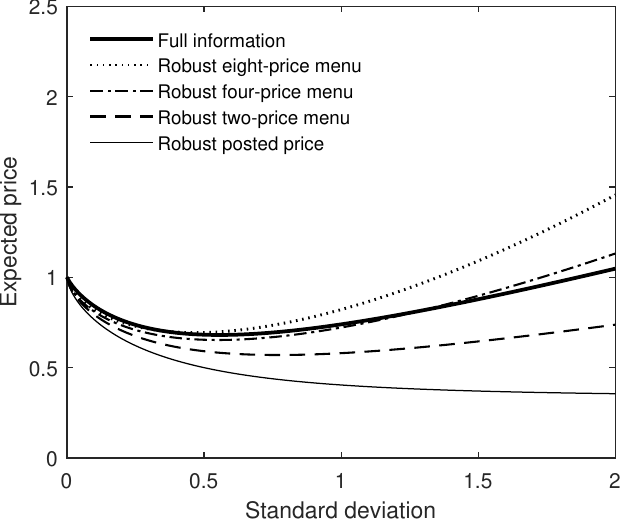}
  \end{subfigure}\hfill
  \begin{subfigure}[b]{0.49\linewidth}
    \centering
    \caption*{Gamma ground-truth distribution}
    \includegraphics[width=\linewidth]{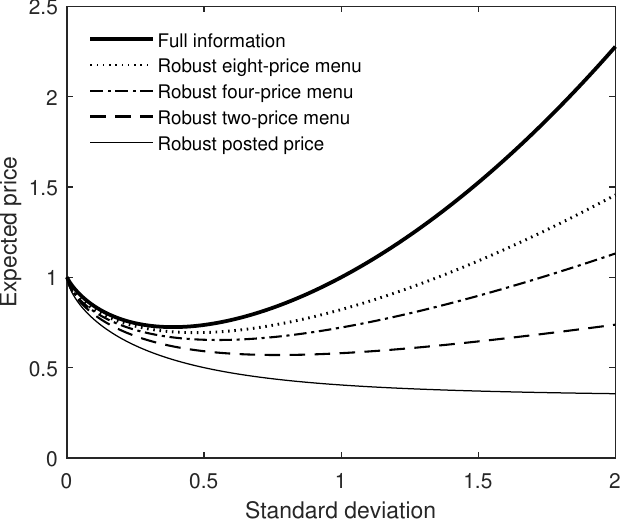}
  \end{subfigure}
  \caption{Expected price under full information and mean--variance partial information, for mean 1.}
  \label{fig:intro}
\end{figure}

Our solution also reveals a pricing phenomenon that a robust posted price cannot capture. To see this, consider first a seller who knows the valuation distribution. Figure~\ref{fig:intro} shows the optimal price of such a seller for lognormal and gamma ground-truth distributions. As the standard deviation increases, the price first decreases and eventually increases above the mean. The same pattern arises for Weibull and beta distributions, among many others. The robust posted price captures only the first part of this response, as it decreases in the standard deviation \citep{chen2022distribution}. As soon as the seller can randomize over two prices, the expected \emph{offered} price also first decreases but eventually increases with the standard deviation. The robust menu thus reproduces this qualitative response of a fully informed seller, while relying only on the mean and the variance. The guaranteed revenue, in contrast, continues to decline, so the upturn reflects a change in the optimal response to heterogeneity rather than a higher revenue guarantee. A single price must compromise between retaining buyers with low valuations and extracting revenue from buyers with high valuations, whereas a menu can address both ends of the market, with a low price that protects against the downside and higher prices that retain access to the upper end.

Our second contribution concerns mean absolute deviation (MAD), another familiar measure of dispersion. The resulting finite-menu problem has a remarkably simple answer. Rather than introducing a second unrelated pricing model, we uncover a structural correspondence: for the relevant worst-case objectives, mean--MAD ambiguity is equivalent to a suitably transformed mean--support problem. This brings Wang's exact geometric-menu theorem into play in a setting where valuations need not have any finite upper bound.
The connection has an unexpected economic consequence. Although MAD admits much heavier upper tails than variance, increasing MAD lowers \emph{every} optimal price. Extreme willingness to pay need not support aggressive pricing when the worst case can make such buyers sufficiently rare. While the mean--variance characterization is the centerpiece of the paper, the MAD result serves two complementary purposes: it identifies a broader equivalence that is useful in its own right, and it makes clearer why the upward price response under variance is not simply a generic benefit of randomizing over several offers.

Our third contribution determines how fast finite menus approach the fully randomized benchmark, and how dispersion affects this speed. Under both mean--variance and mean--MAD information, the revenue gap to full randomization shrinks in inverse proportion to the number of prices, up to a remainder of higher order. While this rate shows that the returns to additional prices diminish, it is the leading coefficient, which we determine exactly, that tells the seller how many prices are needed, as the number of prices required to reach a given share of the fully randomized revenue is approximately proportional to it. We further show that this coefficient increases with dispersion, so that convergence slows down as markets become more volatile. In line with \citet{wang2025power}, a small menu therefore often suffices in calm markets. Diminishing returns do not, however, make a large menu redundant. In volatile markets, the returns to additional prices persist much longer, and the seller may need many prices to come close to the fully randomized revenue. Taken together, the results show that the value of pricing flexibility depends not only on how heterogeneous valuations are, but also on what is known about that heterogeneity.

\subsection{Positioning in the literature}
\paragraph{Robust monopoly pricing with moment information.}
Whereas classical monopoly pricing takes the valuation distribution as given \citep{myerson1981optimal,riley1983optimal}, robust pricing replaces this knowledge with partial information and a worst-case objective, reflecting the appeal of policies that do not rely on an informationally demanding model \citep{Wilson1987}. The distinction between knowing a mean and knowing a mean together with dispersion is consequential: in the unbounded nonnegative model, mean information alone does not secure a positive revenue guarantee \citep{kos2015selling}. Mean--variance \citep{azar2012optimal,chen2022distribution} and mean--MAD information \citep{roos2019chebyshev}, by contrast, yield nontrivial robust posted prices, which decrease as dispersion increases. Adding support information to the mean and variance changes this prescription. \citet{van2024distributionally} show that the optimal posted price then first decreases in the variance and eventually increases, the same pattern we find for the expected offered price. The reasons for the two reversals, however, differ. With a support upper bound, a high variance forces the worst-case distribution to place substantial probability near that bound, which eventually makes a high posted price profitable. Without a support bound, no such force exists, and the upward response appears only once the seller can randomize, already with two prices. The two results therefore complement each other. Starting from mean--variance information, they identify two distinct reasons for a seller to raise prices in more volatile markets, namely a support upper bound and the ability to randomize.

\paragraph{From posted prices to randomized mechanisms.}
A robust seller may benefit from randomizing prices even when a known distribution would make such randomization unnecessary. The unrestricted mechanism-design benchmark with moments and support restrictions has been studied by \citet{carrasco2018optimal}, while \citet{chen2021screening} provide minimax methods for more general convex ambiguity sets. Results for fully randomized pricing under MAD can also be related to \citet{Elmachtoub2021value} and \citet{chen2021screening}. An optimal finite menu, however, cannot be obtained by discretizing the price distributions of these mechanisms, as the seller must choose the prices and their weights jointly, while Nature's worst case depends on the menu itself. Furthermore, the worst-case valuation law can change discontinuously with these choices. This is why a general moment duality argument or a solution for the infinite-menu limit does not by itself settle the finite-menu problem. The converse direction does hold. We show that, as the number of prices grows, the optimal finite menus converge to the fully randomized mechanisms. Table~\ref{tab:overview} positions our results among the existing ones. For each form of information, it lists where the optimal posted price, finite menu, and fully randomized mechanism are characterized. Our results complete the finite-menu row, so that the optimal finite menu is now known for every menu size under all three forms of partial information.
\begin{table}[h!]
\centering
\setlength{\tabcolsep}{4pt}
\begin{tabular}{@{}lccc@{}}
\hline
Mechanism & Mean--variance & Mean--MAD & Mean--support \\ \hline
Posted price & \cite{azar2012optimal} & \cite{roos2019chebyshev} & \cite{kos2015selling} \\
Finite menu & \textbf{Theorem~\ref{th:var_N}} & \textbf{Theorem~\ref{th:mad_N}} & \cite{wang2025power} \\
Full randomization & \cite{carrasco2018optimal} & \cite{chen2021screening} & \cite{carrasco2018optimal} \\ \hline
\end{tabular}
\caption{Overview of existing and new results. The mean--MAD finite-menu characterization is obtained using Theorem \ref{th:bijection} and Wang's mean--support solution.}
\label{tab:overview}
\end{table} 

\paragraph{Alternative objectives and foundations.}
Throughout, our objective is \emph{expected revenue}. The literature also frequently considers the competitive ratio and absolute regret, which compare the expected revenue to that of a fully informed seller, as a ratio and as a difference, respectively. Results under one objective do not automatically carry over to others. Under the competitive ratio, \citet{giannakopoulos2023robust} solve the posted-price problem exactly when only the mean and the variance are known, and give approximations for selling several items. \citet{kleer2024distribution} and \citet{van2025competitive} extend this posted-price analysis to more general dispersion measures, without and with support information, respectively. For randomized mechanisms under the same objective, \citet{wang2024minimax} characterize the optimal fully randomized mechanism under mean--support information and give an approximation under mean--variance information. \citet{wang2025power} then extends these results to finite menus. Absolute-regret formulations appear in \citet{bergemann2008pricing,bergemann2011robust,chen2022distribution,van2025bundling}, while \citet{anunrojwong2024best} compare robustness objectives. More broadly, the moment-based approach connects to distributionally robust optimization, from the seminal inventory model of \citet{scarf1958min} and subsequent work of \citet{gallego1993distribution} to measure-theoretic duality and moment optimization \citep{shapiro2001duality,bertsimas2005optimal,popescu2005semidefinite}.

\subsection{Main findings and organization}
Section~\ref{sec:model} formulates the robust finite-menu pricing problem and introduces the moment-based ambiguity sets. Besides the mean--variance ambiguity set, these include the mean--MAD ambiguity set, and as an auxiliary tool for its analysis, the mean--support ambiguity set.

Section~\ref{sec:variance} presents the optimal finite menus. Theorem~\ref{th:var_N} characterizes robust finite-menu pricing under mean--variance information exactly, for any prescribed number of prices, and is the main result of this paper. It jointly identifies the seller's prices and randomization probabilities together with the structure of the worst-case valuation distribution, so that a problem that initially involves an optimization over all distributions satisfying the moment restrictions ultimately reduces to two scalar equations. The characterization requires neither a bounded valuation support nor equal weights on the menu, and it connects the two benchmarks, recovering the robust posted price for a single price and the fully randomized mechanism as the number of prices grows. Theorem~\ref{th:mad_N} gives the corresponding characterization under mean--MAD information. Its proof relies on a new equivalence between mean--MAD and mean--support ambiguity, which allows the geometric, equally weighted solution of \citet{wang2025power} to be applied in a setting where valuations need not be bounded. The equivalence is new, while the structure of the menu is inherited from the mean--support solution.

Section~\ref{sec:disp} turns to the economic implications. Section~\ref{sec:disp_price} considers how the optimal prices respond to dispersion. With a single posted price, greater variance calls for a lower price. With at least two prices, however, Proposition~\ref{prop:prices_var} shows that the expected price offered to buyers eventually increases with the variance, even as the guaranteed revenue decreases, since the seller can protect the lower end of the market while retaining higher prices for buyers with high valuations. This reversal concerns the expected offered price, rather than the price conditional on a purchase or the seller's revenue, and we distinguish it from the responses of the individual menu prices, which we illustrate numerically. Under mean--MAD information, in contrast, Proposition~\ref{prop:prices_decrease_MAD} shows that all optimal prices decrease as dispersion grows, notwithstanding the possibility of very high valuations, as the worst case can make buyers with extreme valuations too rare to justify targeting them with high prices. The variance finding should therefore not be interpreted as a universal response to more heterogeneous valuations.

Section~\ref{sec:menu_value} then studies the revenue value of additional prices. Building on the importance of menu simplicity highlighted by \citet{wang2025power}, we ask not only how many prices are useful, but also how the information about valuations shapes the gains from expanding the menu. Theorem~\ref{th:convergence_rate} determines how fast the optimal revenue approaches that of full randomization under both forms of information and how this speed depends on dispersion, which yields a simple rule for the number of prices needed to secure a given share of the fully randomized revenue.

Section~\ref{sec:variance_proof} develops the full proof of Theorem~\ref{th:var_N}, which forms the technical core of the paper and makes each step of the sketch in Section~\ref{sec:variance} rigorous. Section~\ref{sec:bijection} proves Theorem~\ref{th:mad_N} through an equivalence between mean--MAD and mean--support ambiguity, which we establish in Theorem~\ref{th:bijection} and prove in Section~\ref{sec:proof}. For a broad class of objectives that includes the revenue of any menu, the worst case under mean--MAD ambiguity coincides with that under a suitably chosen mean--support ambiguity set, so that Theorem~\ref{th:mad_N} follows from the known solution of the mean--support pricing problem. Section~\ref{sec:conclusion} concludes with the implications of our results for the choice of dispersion measure and with directions for future research. The appendices contain the remaining proofs, including the convergence of the optimal finite menus to the fully randomized mechanisms.

\section{Problem formulation}\label{sec:model}
We consider a monopolist selling a single item to a buyer with a non-negative random valuation $X$ distributed according to an unknown distribution $\bP$. The seller may randomize over $N \in \N$ prices. A feasible menu $(w,p)$ consists of prices $0<p_1\leq\dots\leq p_N$ and weights in
\[
    \Delta_N=\Bigl\{w\in[0,1]^N:\sum_{i=1}^Nw_i=1\Bigr\},
\]
and we write
\[
    \mathcal H=\Bigl\{(w,p)\in\Delta_N\times\R_+^N:0<p_1\leq\dots\leq p_N\Bigr\}
\]
for the set of feasible menus. We then aim to maximize expected revenue, i.e.,
\begin{align}\label{eq:rev}
    \sum_{i=1}^N w_i p_i \bP(X \ge p_i).
\end{align}
Since the seller does not know $\bP$ exactly, we model this partial information by considering $\bP \in \cP$, where $\cP$ is an ambiguity set of possible distributions. The seller then chooses the menu to maximize revenue in the worst case, where the distribution is selected from $\cP$ by an adversary that we refer to as Nature. For fixed prices $0<p_1\leq\dots\leq p_N$, the revenue the seller can guarantee by choosing the weights optimally is
\begin{align}\label{eq:Rp}
    R(p)=\sup_{w\in\Delta_N}\ \inf_{\bP\in\cP}\ \sum_{i=1}^Nw_ip_i\bP(X\geq p_i),
\end{align}
and optimizing over the prices as well gives the distributionally robust pricing problem
\begin{align}\label{eq:Npricing_mad}
    R^*=\sup_{p}R(p)=\sup_{(w,p)\in\mathcal{H}}\inf_{\bP\in\cP}\sum_{i=1}^Nw_ip_i\bP(X\geq p_i).
\end{align}
We sometimes write $R^*(N)$ instead of $R^*$ when the dependence on the menu size matters.

Our main analysis imposes the two moments most commonly used to summarize a market: the mean and variance of willingness to pay. These moments do not identify a demand curve, and their interaction matters. In particular, a fixed mean without further information allows Nature to assign arbitrarily small probability to arbitrarily high valuations, which precludes a positive revenue guarantee even when the seller can randomize over a continuum of prices \citep{kos2015selling}. Prescribing the second moment limits this possibility while leaving the seller's optimal finite-menu response far from obvious. We also consider MAD information, which we analyze through its equivalence with mean--support information.

We begin with variance, which penalizes large deviations quadratically and defines the principal ambiguity set
$$
\cP_\sigma(\mu,\sigma) = \Big\{\bP:
\bP(X \in [0,\infty))=1,\;
\bE_{\bP}[X]=\mu,\;
\bE_{\bP}[X^2]= \mu^2+\sigma^2\Big\},
$$
where $\mu>0$ denotes the mean and $\sigma>0$ the standard deviation.

As a complementary benchmark, mean absolute deviation (MAD) penalizes deviations linearly and may remain finite even when variance is infinite. Its ambiguity set is
$$
\cP_d(\mu,d) = \Big\{\bP:
\bP(X \in [0,\infty))=1,\;
\bE_{\bP}[X]=\mu,\;
\bE_{\bP}[|X-\mu|]= d\Big\},
$$
where $\mu>0$ denotes the mean and $d\in(0,2\mu)$ the MAD. This upper bound on the MAD is not an additional assumption. Since $\bE[X-\mu]=0$, the MAD satisfies $d=\bE[|X-\mu|]=2\,\bE[(\mu-X)^+]$. For nonnegative valuations, $(\mu-X)^+\leq\mu$ with equality if and only if $X=0$, which cannot hold almost surely since $\bE[X]=\mu>0$. Consequently, $d<2\mu$.

Our analysis of MAD information relies on the mean--support ambiguity set, in which the MAD constraint is replaced by a known finite upper bound $\beta$. It is defined as
$$
\cP_\beta(\mu,\beta) = \Big\{\bP:
\bP(X \in [0,\beta])=1,\;
\bE_{\bP}[X]=\mu\Big\},
$$
where $0<\mu<\beta < \infty$. All three ambiguity sets fix the mean and limit how much probability Nature can place on high valuations, but do so in different ways. These differences shape the worst case, and with it the optimal menu.

\section{Optimal finite menus}\label{sec:variance}
As our main contribution, Theorem~\ref{th:var_N} solves problem~\eqref{eq:Npricing_mad} exactly under mean--variance ambiguity for every menu size, yielding the optimal prices, their weights, and the optimal expected revenue. The entire menu is built from its two lowest prices, which follow from a system of two equations in two unknowns.

\begin{theorem}\label{th:var_N}
Consider $\cP = \cP_\sigma(\mu,\sigma)$ with $N\geq2$. For a pair $(p_1,p_2)$, define
\begin{align}\label{eq:zdef}
    z=\frac12\left(p_1+\frac{\mu^2+\sigma^2}{p_1}\right)
\end{align}
and the prices
\begin{align}\label{eq:rec}
    p_{i}=z-\sqrt{(z-p_{i-1})(z+p_{i-1}-2r_{i-2}p_{i-1})},\qquad i=3,\dots,N,
\end{align}
where $r_i=p_{i+1}/p_{i}$ for $i=1,\dots,N$ with the convention that $p_{N+1}=z$. The system of two equations
\begin{align}
    z&=\bigl(2r_{N-1}-1\bigr)p_N,\label{eq:terminal}\\
    \bigl(\mu^2+\sigma^2-p_1^2\bigr)\sum_{i=1}^N\bigl(r_i-1\bigr)&=\bigl(\mu-p_1\bigr)\sum_{i=1}^Np_i\bigl(r_i^2-1\bigr),\label{eq:moment}
\end{align}
in the two unknowns $p_1$ and $p_2$ has a unique solution $(p^*_1,p^*_2)$ for which $0<p^*_1<\dots<p^*_N$ with $p^*_3,\dots,p^*_N, z^*,r^*_1,\dots,r^*_N$ the corresponding quantities. Then the prices $p^*_1,\dots,p^*_N$ with weights
\begin{align}\label{eq:var_menu}
    w^*_i=\frac{\bigl(r^*_i-1\bigr)\bigl(2z^*-p^*_i-p^*_{i+1}\bigr)}{\sum_{j=1}^N\bigl(r^*_j-1\bigr)\bigl(2z^*-p^*_j-p^*_{j+1}\bigr)},\qquad i=1,\dots,N,
\end{align}
form an optimal menu with expected revenue
\begin{align}\label{eq:rev_var_N}
    R^*=\frac{\mu-p^*_1}{\sum_{i=1}^N\bigl(r^*_i-1\bigr)}.
\end{align}
\end{theorem}

Theorem~\ref{th:var_N} connects the two extremes of the pricing problem. For $N=1$, the characterization reduces to finding $p^*_1$ and condition~\eqref{eq:terminal} dropping out, which recovers the robust posted price of \citet{azar2012optimal}. As $N\to\infty$, the optimal menu converges to that of the fully randomized mechanism of \citet{carrasco2018optimal}, which we prove in Appendix~\ref{ap:asympt}.

We now sketch the proof of Theorem~\ref{th:var_N}, which shows where its conditions come from and how the difficulties of the mean--variance case are resolved. The full proof is presented in Section~\ref{sec:variance_proof}. Problem~\eqref{eq:Npricing_mad} under mean--variance ambiguity is a zero-sum game between the seller, who chooses a menu $(w,p)\in\mathcal H$ to maximize expected revenue, and an adversarial Nature, who chooses a valuation distribution $\bP\in\cP_\sigma(\mu,\sigma)$ to minimize it. Fix the prices $p$, and consider the game in which the seller only chooses the weights. Its value is the revenue $R(p)$ that the seller can guarantee with these prices, and it has a saddle point, i.e., a Nash equilibrium in which the weights are a best response to Nature's distribution and vice versa \citep[see][for the same view under full randomization]{carrasco2018optimal}. Nature's equilibrium distribution is only approached in the limit, as we discuss below. Every price receives positive weight, and consequently, the seller's best response requires that every price on the menu earns the same revenue,
\begin{align}\label{eq:equal_revenue}
    p_i\,\bP(X\geq p_i) = R(p),\qquad i=1,\dots,N,
\end{align}
as the seller would otherwise shift weight towards the price that earns the most. Nature's equilibrium distribution thus leaves the seller indifferent between the prices, and the probability of a sale at each price $p_i$ equals $R(p)/p_i$.

The equal-revenue property determines the worst case almost completely. It fixes how much probability Nature places between consecutive prices. A valuation just below a price contributes as much as possible to the mean without leading to a sale at that price, so Nature pushes the probability in each interval up against the next price. The remaining probability, that of a buyer who purchases at every price on the menu, is concentrated on a single high valuation $z(p)$ above all prices. Pushing the probability in each interval ever closer to the next price, Nature approaches its worst case through distributions that converge to
\begin{align}\label{eq:worstcase_sketch}
    \bP^*=\Bigl(1-\frac{R(p)}{p_1}\Bigr)\delta_{p_1}+\sum_{i=2}^NR(p)\Bigl(\frac1{p_{i-1}}-\frac1{p_i}\Bigr)\delta_{p_i}+\frac{R(p)}{p_N}\,\delta_{z(p)},
\end{align}
where $\delta_x$ denotes the point mass at $x$.

The high valuation $z(p)$ is the main source of difficulty. The mean pushes probability towards high valuations while the variance penalizes it, and Nature settles on a valuation $z(p)$ that is neither a price nor fixed by the information, but determined jointly with the menu. Consequently, $\bP^*$ is not yet fully known, as it still depends on $z(p)$ as well as on the revenue $R(p)$. Since $\bP^*$ describes Nature's worst case only if it has mean $\mu$ and variance $\sigma^2$, these two moment constraints pin down $R(p)$ and $z(p)$ jointly. With $p_{N+1}=z(p)$ and $r_i=p_{i+1}/p_i$, they read
\begin{align}\label{eq:two_moments}
    \mu = p_1+R(p)\sum_{i=1}^N(r_i-1),
    \qquad
    \mu^2+\sigma^2 = p_1^2+R(p)\sum_{i=1}^Np_i(r_i^2-1),
\end{align}
and we refer to $\bP^*$ with $R(p)$ and $z(p)$ satisfying \eqref{eq:two_moments} as the limiting worst-case distribution. At the optimal prices, solving the first equation for $R(p)$ gives \eqref{eq:rev_var_N}, and inserting \eqref{eq:rev_var_N} into the second equation gives \eqref{eq:moment}.

It remains to optimize over the prices, which is the main technical step of the proof. Under mean--support information, a similar analysis gives the worst-case revenue as an explicit function of the prices without any implicitly determined component. Its maximization is then straightforward, as the first-order conditions can be solved directly \citep{wang2025power}. Under mean--variance information, the worst-case revenue also depends on $z(p)$, which is only determined implicitly by \eqref{eq:two_moments}. Its maximization is therefore considerably more difficult, and in general a non-convex problem. By \eqref{eq:two_moments}, the worst-case revenue satisfies $R(p)=\Psi(p,z(p))$ for an explicit function $\Psi$ of the prices and $z$. Rather than solving \eqref{eq:two_moments} for $z(p)$, we treat $z$ as an additional variable, so that maximizing $\Psi$ over both the prices and $z$ gives an upper bound on $R^*$. This maximization over $N+1$ variables is not straightforward. Its domain is not compact, since $p_1$ may approach zero and $z$ may grow without bound, and $\Psi$ is not concave, so that a maximum need not exist, and even if it does, the first-order conditions need not identify it. We therefore proceed in three steps. First, we show that the maximum is attained, since at points $(p,z)$ where $\Psi$ is close to its supremum, the lowest price is bounded away from zero and $z$ is bounded from above. Second, we show that every maximizer is an interior point, i.e., has $p_1<p_2<\dots<p_N<z$, so that it satisfies the first-order conditions, which give \eqref{eq:zdef}, \eqref{eq:rec} and \eqref{eq:terminal}. Third, a one-dimensional monotonicity argument shows that these conditions have at most one solution. Together, the three steps imply that the unique solution of the first-order conditions is the global maximizer.

This maximum, however, is only an upper bound on $R^*$, as it treats $z$ as a free variable. If the maximizer $(p^*,z^*)$ satisfied $z^*=z(p^*)$, the bound would be attained by the prices $p^*$, but since we did not solve for $z(p)$, this cannot be verified directly. We therefore derive a matching lower bound. If the upper bound were tight, the analysis for fixed prices would determine the optimal weights for the prices $p^*$, which are given in \eqref{eq:var_menu}. Together with the prices $p^*$, they define a candidate menu, whose guaranteed revenue is a lower bound on $R^*$. By constructing a feasible solution of the dual of Nature's problem, we show that this candidate menu guarantees exactly the maximum of $\Psi$. The lower and upper bounds therefore coincide, which confirms that the bound is tight and that the candidate menu is indeed optimal.

We now turn to mean--MAD ambiguity, whose solution we use alongside Theorem~\ref{th:var_N} in the remainder of the paper.

\begin{theorem}\label{th:mad_N}
Consider problem~\eqref{eq:Npricing_mad} with $\cP = \cP_d(\mu,d)$ and $N\in\N$. An optimal menu $(w^*,p^*)$ is given by
\begin{align}\label{eq:mad_menu}
    w^*_i = \frac1N,\qquad p^*_i = p^*_1\left(\frac{\mu}{p^*_1}\right)^{\frac{i-1}{N}},\qquad i=1,\dots,N,
\end{align}
where $p^*_1$ is the unique solution on $(0,\mu-d/2)$ to
\begin{align}\label{eq:p1}
    \frac{\mu-d/2}{p^*_1}+N\left(\frac{p^*_1}{\mu}\right)^{1/N}=N+1,
\end{align}
and its expected revenue is
\begin{align}\label{eq:rev_N_2a}
    R^* = \frac{p^*_1p^*_N}{\mu} = \frac{N+1}{N}\,p^*_1-\frac{\mu-d/2}{N}.
\end{align}
\end{theorem}

Compared with Theorem~\ref{th:var_N}, the optimal menu is considerably simpler. The seller randomizes uniformly over geometrically spaced prices and the entire menu follows from a single equation for its lowest price, and the optimal expected revenue is the product of the lowest and the highest price, divided by the mean. Theorem~\ref{th:mad_N} connects the two extremes of the pricing problem as well. For $N=1$, it recovers the robust posted price of \citet{roos2019chebyshev}, and as $N\to\infty$, the optimal menu converges to the fully randomized mechanism of \citet{chen2021screening}, which we prove in Appendix~\ref{ap:asympt}. The proof of Theorem~\ref{th:mad_N} rests on an equivalence between mean--MAD and mean--support ambiguity, which we present in Section~\ref{sec:bijection}.

\section{Pricing implications of valuation dispersion}\label{sec:disp}
The structural solution to mean--variance finite-menu pricing is interesting for more than its tractability. It reveals that the availability of even one additional price changes how a seller should respond to greater valuation heterogeneity. The usual posted-price intuition emphasizes downside protection and therefore favors a lower price. Randomization creates a second margin: the seller can keep a low offer while assigning probability to an offer aimed at buyers with higher willingness to pay. We begin with this mean--variance phenomenon, carefully separating the behavior of the expected offered price from guaranteed revenue. We then examine the optimal menu itself and, only after establishing the main mechanism, return to MAD as an illuminating counterexample to a universal upward-pricing effect. This extends Wang's insight about the performance of simple menus to a different issue: their effect on the direction of optimal pricing decisions.

\subsection{Price responses to valuation dispersion}\label{sec:disp_price}
\begin{figure}[h!]
  \centering
  \begin{subfigure}[b]{0.48\linewidth}
    \centering
    \includegraphics[width=\linewidth]{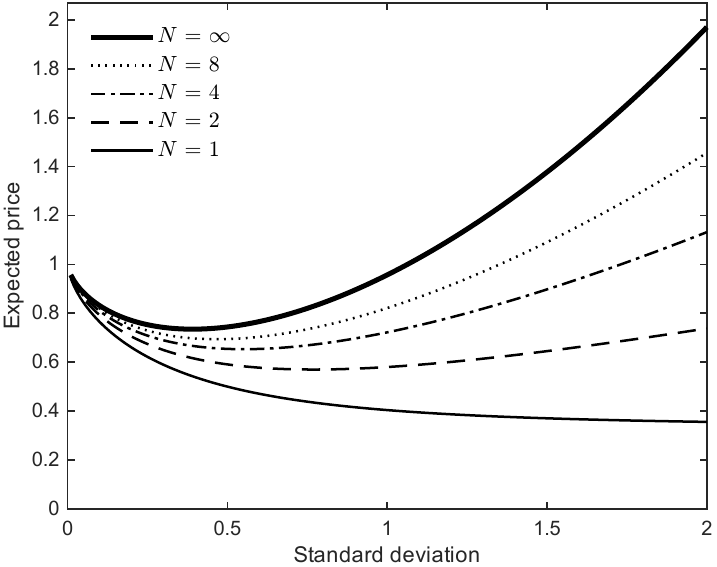}
  \end{subfigure}\hfill
  \begin{subfigure}[b]{0.48\linewidth}
    \centering
    \includegraphics[width=\linewidth]{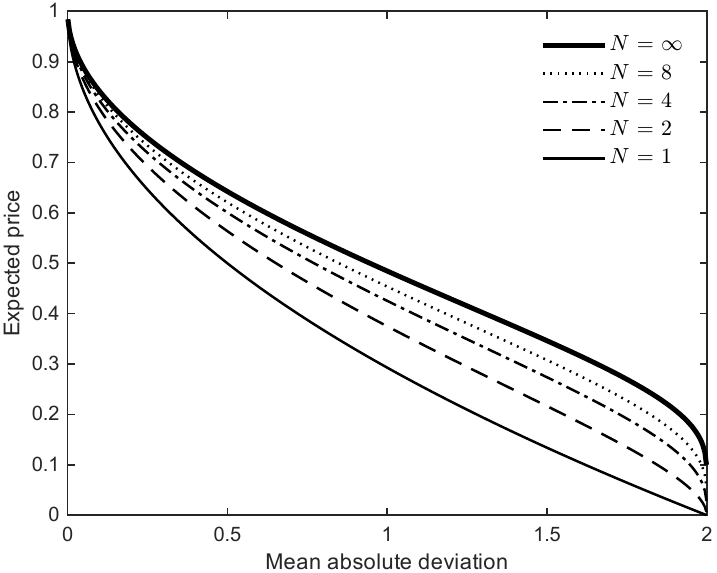}
  \end{subfigure}
    \caption{Expected price of the optimal menu for different values of $N$, for $\cP=\cP_\sigma(\mu,\sigma)$ as a function of $\sigma$ (left) and for $\cP=\cP_d(\mu,d)$ as a function of $d$ (right), with $\mu=1$.}
  \label{fig:expected_prices}
\end{figure}

Figure~\ref{fig:expected_prices} displays the main mean--variance result on the left and plots the expected \emph{offered} price, $\bE(P;\sigma)=\sum_{i=1}^Nw^*_ip^*_i$. Under variance information, the initial decline gives way to an increase when the seller can offer multiple prices. The MAD comparison, shown alongside for perspective, behaves differently. The two panels use different dispersion statistics, $\sigma$ and $d$; their horizontal scales should not be compared directly.

For a single posted price, the two measures agree qualitatively: the optimal price decreases as the standard deviation increases \citep{chen2022distribution}, and likewise as MAD increases \citep{roos2019chebyshev}. A single price must balance the risk of losing buyers with low valuations against the revenue from buyers willing to pay more. A menu can separate these responses. We track both its individual prices and the expected offered price $\bE(P;\sigma)$, which reduces to the posted price for $N=1$. A higher offered price is not, by itself, a higher revenue guarantee.
 
\begin{figure}[h!]
  \centering
  \begin{subfigure}[b]{0.48\linewidth}
    \centering
    \caption*{$N = 2$}
    \includegraphics[width=\linewidth]{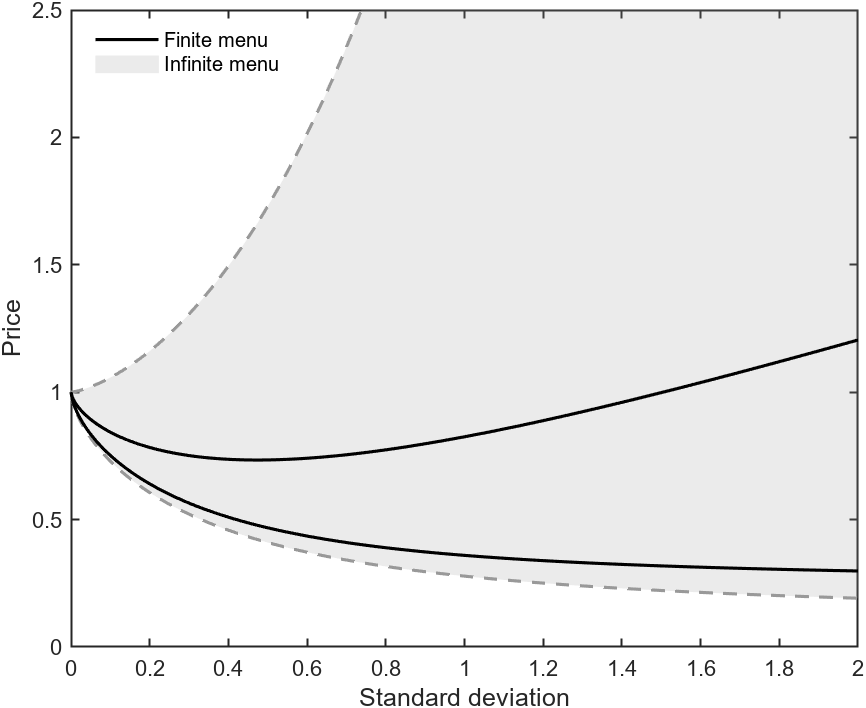}
  \end{subfigure}
  \begin{subfigure}[b]{0.48\linewidth}
    \centering
    \caption*{$N = 4$}
    \includegraphics[width=\linewidth]{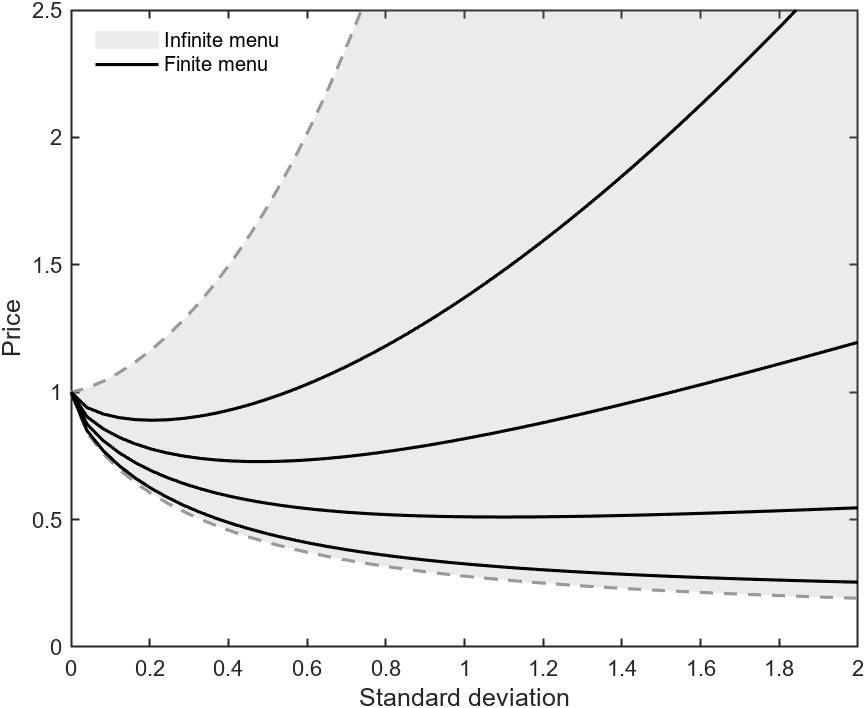}
  \end{subfigure}\hfill
  \caption{Optimal finite-menu prices for $\cP = \cP_\sigma(\mu,\sigma)$ and different values of $N$. The finite-menu prices are compared to the infinite-menu prices, as functions of the standard deviation with $\mu=1$.}
  \label{fig:prices_var}
\end{figure}

The mean--variance case, our primary focus, is different from the posted-price benchmark. Figure~\ref{fig:prices_var} shows the optimal menu widening as dispersion increases: numerically, the lowest price moves down, while higher prices turn upward and can exceed the mean. This suggests a two-sided response, retaining a low-price option while using another price to reach buyers at the upper end of the worst-case market. The next proposition, proved in Appendix~\ref{ap:prices_var}, makes the robust conclusion precise for the expected offered price and worst-case revenue. The stronger apparent pattern for each individual high price is illustrated numerically, not claimed as a theorem.
 
\begin{proposition}\label{prop:prices_var}
    Let $N \geq 2$ and $\cP = \cP_\sigma(\mu,\sigma)$.
    \begin{enumerate}
        \item[(i)] There are thresholds $0<\sigma_l\leq\sigma_h<\infty$
        such that the expected price $\bE(P;\sigma)$ is decreasing in $\sigma$ on
        $(0,\sigma_l)$ and increasing in $\sigma$ on $(\sigma_h,\infty)$.
        \item[(ii)] The optimal expected revenue $R^*$ is decreasing in $\sigma$.
    \end{enumerate}
\end{proposition}

Proposition~\ref{prop:prices_var} establishes the price reversal without asserting its exact location. Between $\sigma_l$ and $\sigma_h$, the behavior of the expected offered price is not resolved analytically. The reason is that the optimal prices depend on $\sigma$ only implicitly through \eqref{eq:terminal}--\eqref{eq:moment}, and are linked to each other through the recursion \eqref{eq:rec}.  Across the numerical instances examined, the price appears to turn exactly once, suggesting the conjecture $\sigma_l=\sigma_h$ for every $N\geq2$. This sharper single-turning-point property remains open.
 
Mean--MAD provides an illuminating benchmark for this variance result. In that setting, the posted-price intuition continues to hold for the entire menu. The next proposition, proved in Appendix~\ref{ap:monotonicity}, shows that every optimal price decreases as the MAD increases, for every menu size, and that the guaranteed revenue decreases as well. Figure~\ref{fig:prices_mad} in Appendix~\ref{ap:MADdecrease} illustrates this uniform downward movement of the menu.

\begin{proposition}\label{prop:prices_decrease_MAD}
    Let $N \in \N$ and $\cP = \cP_d(\mu,d)$.
    \begin{enumerate}
        \item[(i)] The optimal prices $p^*_1,\dots,p^*_N$, and hence the expected price $\bE(P;d)$, are decreasing in $d$.
        \item[(ii)] The optimal expected revenue $R^*$ is decreasing in $d$.
    \end{enumerate}
\end{proposition}

The geometric menu and equal weights are not, in themselves, a departure from \citet{wang2025power}: they follow from the mean--support solution once we establish the mean--MAD/mean--support equivalence in Theorem~\ref{th:bijection}. What is distinctive here is that this inherited menu structure produces a strong comparative-statics conclusion under a dispersion measure that allows arbitrarily high valuations. More prices improve the guarantee, but do not reverse the downward price response.
 
The economic significance of the variance result is not that additional variance improves the revenue guarantee: it does not. Rather, limited flexibility allows the seller to manage downside exposure and high-price opportunities within the same menu, overturning the posted-price comparative statics for the expected offered price. This sharpens the connection to \citet{wang2025power}, who establishes the substantial \emph{value} of small menus; our mean--variance characterization shows that introducing a second price can also alter the \emph{direction} in which pricing responds to heterogeneity. The comparison with MAD further shows that this reversal depends on the structure of the moment restriction, not only on the availability of randomization.
 
The variance mechanism follows from the structure of Nature's response. Theorem~\ref{th:var_N} identifies a discrete optimal worst case with a finite high support point $z(p)$ above the menu, with the remaining support tied to the seller's prices. The second-moment condition limits the ease with which Nature can move mass arbitrarily far into the upper tail. The seller's menu can therefore combine a lower price that protects the downside with higher prices responding to the upper part of this worst case. Unlike Wang's mean--support menu, the optimal variance menu is characterized by a nonlinear recursion and generally nonuniform weights; the reversal in the expected offered price is a consequence of this different structure. This comparison is about optimal worst-case mechanisms under distinct information sets, not about a literal nesting of variance and MAD ambiguity sets.

Why does MAD behave differently? Under MAD, Nature may satisfy the dispersion constraint through exceptionally high valuations that occur with exceptionally small probabilities. The presence of an extreme upper tail therefore need not give the seller a profitable high-price target: the seller risks pricing itself out of the market. Theorem~\ref{th:bijection} makes this intuition exact for menu revenue. Despite the possibility of unbounded valuations, the mean--MAD problem is equivalent to a mean--support problem with upper bound $\mu$ and effective mean $\mu-d/2$. Increasing MAD lowers this effective mean without increasing the effective upper bound. The resulting optimal menu is precisely Wang's geometric, equal-weight menu at the transformed parameters.

Finally, the two settings agree on an important point established more broadly by \citet{wang2025power}: a limited number of prices can already capture a substantial share of the gains from randomization. Our additional question is how dispersion changes that share, and hence the menu size needed for a given revenue guarantee. Section~\ref{sec:disp_menu} studies this issue quantitatively.

\subsection{The value of additional prices}\label{sec:menu_value}\label{sec:disp_menu}
A posted price is the simplest mechanism available to the seller, but it also earns the least. At the other extreme, full randomization earns the most, but it is also the most complex, as it requires the seller to randomize over a continuum of prices. This leads us to investigate the marginal value of finite menus.

Under mean--variance ambiguity, the optimal expected revenue of the fully randomized mechanism \citep{carrasco2018optimal} is the unique solution $\pi_\sigma\in(0,\mu)$ of
\begin{align}\label{eq:pisigma}
    \mu^2+\sigma^2 = \pi_\sigma^2\bigl(2e^{\mu/\pi_\sigma-1}-1\bigr),
\end{align}
while under mean--MAD ambiguity \citep{Elmachtoub2021value,chen2021screening}, it is the unique solution $\pi_d\in(0,\mu-d/2)$ of
\begin{align}\label{eq:pid}
    \mu-\frac d2 = \pi_d\Bigl(1+\ln\frac{\mu}{\pi_d}\Bigr).
\end{align}
As $N$ grows large, $R^*(N)$ converges to $\pi_\sigma$ and $\pi_d$, respectively. Appendix~\ref{ap:asympt} collects the corresponding asymptotic statements and proofs. Table~\ref{tab:rev} then reports the ratios $R^*(N)/\pi_\sigma$ and $R^*(N)/\pi_d$ for various values of $N$, $\sigma$ and $d$, showing how rapidly finite menus approach the fully randomized benchmark.

Two patterns emerge under both dispersion measures. First, for a fixed level of dispersion, the gains from enlarging the menu diminish with $N$. Second, convergence to the optimal expected revenue becomes slower when dispersion is higher. Taken together, these observations show that a small menu often suffices, since the gains from adding more prices diminish quickly, but when dispersion is high, a larger menu may nevertheless be required to capture a substantial share of the optimal expected revenue. The next theorem makes both patterns precise.

\begin{theorem}[Rate of convergence]\label{th:convergence_rate}\leavevmode\begin{enumerate}
        \item[(i)] Let $\cP=\cP_\sigma(\mu,\sigma)$, and let $\pi_\sigma$ be as in \eqref{eq:pisigma}. The convergence of $R^*(N)$ to $\pi_\sigma$ is linear in $1/N$, i.e.,
        $$R^*(N)=\pi_\sigma\left(1-\frac{\kappa_\sigma}{N}+O\left(\frac{\ln(N)}{N^2}\right)\right), \qquad \kappa_\sigma=\frac{2\Bigl(\frac12\ln\frac{1+\sqrt{a_\sigma}}{1-\sqrt{a_\sigma}}-\sqrt{a_\sigma}\Bigr)^2}{\mu/\pi_\sigma-1-a_\sigma},$$
        with $a_\sigma=1-e^{1-\mu/\pi_\sigma}$. Moreover, the coefficient $\kappa_\sigma$ is increasing in $\sigma$.
        \item[(ii)] Let $\cP=\cP_d(\mu,d)$, and let $\pi_d$ be as in \eqref{eq:pid}. The convergence of $R^*(N)$ to $\pi_d$ is linear in $1/N$, i.e.,
        $$R^*(N)=\pi_d\left(1-\frac{\kappa_d}{N}+O\left(\frac{1}{N^2}\right)\right), \qquad \kappa_d=\frac{\ln(\mu/\pi_d)}{2}.$$
        Moreover, the coefficient $\kappa_d$ is increasing in $d$.
    \end{enumerate}
\end{theorem}

The proof is found in Appendix~\ref{ap:convergence_rate}. Theorem~\ref{th:convergence_rate} thus confirms the patterns observed in Table~\ref{tab:rev}. The $1/N$ rate captures the diminishing returns from enlarging the menu, and the leading coefficient increasing in the dispersion parameter captures the need for a larger menu at higher dispersion levels. 

The coefficient also tells the seller how many prices the menu should contain. Let $N_{1-\varepsilon}$ denote the smallest menu size for which $R^*(N_{1-\varepsilon})$ captures at least $1-\varepsilon$ of the optimal expected revenue. By Theorem~\ref{th:convergence_rate}, this equals $1-\kappa_\sigma/N$ under mean--variance ambiguity and $1-\kappa_d/N$ under mean--MAD ambiguity, up to a remainder of smaller order in $1/N$. Ignoring this remainder, the fraction is at least $1-\varepsilon$ precisely when $N\geq\kappa_\sigma/\varepsilon$ or $N\geq\kappa_d/\varepsilon$, which gives the approximations $N_{1-\varepsilon}\approx\kappa_\sigma/\varepsilon$ and $N_{1-\varepsilon}\approx\kappa_d/\varepsilon$. The last row of each panel in Table~\ref{tab:rev} reports the exact value of $N_{0.95}$, computed from the optimal menus. For the dispersion levels in Table~\ref{tab:rev}, the approximations $\kappa_\sigma/0.05$ and $\kappa_d/0.05$ deviate from these exact values by less than two prices.

\begin{table}[h!]
\centering
\begin{tabular}{c cccccc c cccccc}
\toprule
 & \multicolumn{6}{c}{Mean--variance, $\sigma$} & & \multicolumn{6}{c}{Mean--MAD, $d$} \\
\cmidrule(lr){2-7} \cmidrule(lr){9-14}
$N$ & 0.2 & 0.5 & 1.0 & 1.5 & 2.0 & 4.0 & & 0.2 & 0.5 & 1.0 & 1.5 & 1.8 & 1.9 \\
\midrule
1  & 0.81 & 0.61 & 0.38 & 0.25 & 0.18 & 0.06 & & 0.80 & 0.65 & 0.46 & 0.27 & 0.13 & 0.07 \\
2  & 0.89 & 0.76 & 0.59 & 0.48 & 0.40 & 0.25 & & 0.89 & 0.80 & 0.67 & 0.51 & 0.37 & 0.28 \\
4  & 0.94 & 0.86 & 0.76 & 0.68 & 0.62 & 0.49 & & 0.94 & 0.89 & 0.81 & 0.72 & 0.61 & 0.54 \\
8  & 0.97 & 0.93 & 0.87 & 0.82 & 0.78 & 0.70 & & 0.97 & 0.94 & 0.90 & 0.85 & 0.78 & 0.74 \\
16 & 0.98 & 0.96 & 0.93 & 0.90 & 0.88 & 0.84 & & 0.98 & 0.97 & 0.95 & 0.92 & 0.89 & 0.86 \\
32 & 0.99 & 0.98 & 0.96 & 0.95 & 0.94 & 0.91 & & 0.99 & 0.99 & 0.97 & 0.96 & 0.94 & 0.93 \\
\midrule
$N_{0.95}$ & 6 & 13 & 24 & 33 & 39 & 57 & & 5 & 10 & 17 & 27 & 38 & 47 \\
\bottomrule
\end{tabular}
\caption{$R^*(N)/\pi_\sigma$ and $R^*(N)/\pi_d$ for various $N$, $\sigma$, and $d$ with $\mu = 1$. $N_{0.95}$ reports the smallest $N$ for which the ratio is at least $0.95$.}
\label{tab:rev}
\end{table}

\section{Proof of the Mean--Variance Finite Menu Theorem}\label{sec:variance_proof}
Throughout this proof, with a slight abuse of notation, we write $R(p)$ and $R^*$ for the values in \eqref{eq:Rp} and \eqref{eq:Npricing_mad}, respectively, under the strict-tail convention $\bP(X>p_i)$ rather than the weak-tail convention $\bP(X\geq p_i)$. The following lemma, proven in Appendix~\ref{ap:Lemma_taileq}, shows that this leaves $R^*$ unchanged, and that optimal menus remain optimal.

\begin{lemma}\label{lem:tail_equivalence}
Problem \eqref{eq:Npricing_mad} has the same value under the strict-tail convention, and every menu that is optimal under the strict-tail convention is optimal for \eqref{eq:Npricing_mad}.
\end{lemma}

The remainder of the proof consists of four steps. First, we fix the prices and show, using Sion's minimax theorem, that Nature's worst-case distribution places masses on the intervals $[0,p_1]$, $(p_i,p_{i+1}]$ and $(p_N,\infty)$ that are determined by the prices up to the single unknown $R(p)$. Second, we use duality to show that the mass in each of these intervals concentrates at its right endpoint, and for $(p_N,\infty)$ at a point $z(p)$. The mean and the variance constraints then give $R(p)=\Psi(p,z(p))$ for an explicit function $\Psi$. Since $z(p)$ is known only implicitly, we treat $z$ as a free variable, which results in the upper bound $R^*\leq\sup_{p,z}\Psi(p,z) = \Psi^*$. Third, we show that $\Psi^*$ is attained at a unique point $(p^*,z^*)$, which is characterized by the conditions of Theorem~\ref{th:var_N}. Finally, we derive candidate weights $w$ and show, by constructing a feasible dual solution, that the menu $(p^*,w)$ guarantees $\Psi^*$, so that $\Psi^* = R^*$.

Fix prices $0<p_1\leq\dots\leq p_N$ with $R(p)>0$. This is without loss of generality, since $R^*>0$. Indeed, for $0<t<\mu$, Cantelli's inequality gives $\bP(X>\mu-t)\ge t^2/(t^2+\sigma^2)$ for every $\bP\in\cP_\sigma(\mu,\sigma)$, so that the menu charging $p_1 = \mu-t$ with $w_1 = 1$ guarantees a positive revenue. The set of revenue vectors that Nature can induce,
\[
    \mathcal A(p):=\bigl\{(p_i\bP(X>p_i))_{i=1}^N:\bP\in\cP_\sigma(\mu,\sigma)\bigr\},
\]
is bounded and convex, so that its closure $\bar{\mathcal A}(p)$ is compact and convex. Consequently, by the minimax theorem \citep{sion1958general} we obtain
\begin{align}\label{eq:minimax}
    R(p)=\max_{w\in\Delta_N}\min_{a\in\bar{\mathcal A}(p)}\sum_{i=1}^Nw_ia_i=\min_{a\in\bar{\mathcal A}(p)}\max_{w\in\Delta_N}\sum_{i=1}^Nw_ia_i
\end{align}
together with a saddle point $(w^*,a^*)$. 

We may assume without loss of generality that $w_i^*>0$ for all $i=1,\dots,N$. Otherwise, replace every price $p_i$ with $w^*_i=0$ by the largest price $p_j<p_i$ with $w^*_j>0$, or by the smallest such $p_j>p_i$ if there is none, and divide the weight $w^*_j$ equally among its copies. Denote the new prices and weights by $\bar p$ and $\bar w$, so that $0<\bar p_1\le\dots\le\bar p_N$ and $\bar w_i>0$ for all $i$. Since each $\bar p_i$ is one of the prices $p_j$, every menu with prices $\bar p$ randomizes over a subset of the prices in $p$, so that $R(\bar p)\le R(p)$. Conversely, let $C_j$ consist of $j$ and the indices $i$ whose price $p_i$ was replaced by $p_j$, so that $\bar w_i=w_j^*/|C_j|$ for $i\in C_j$. Then, for every $\bP\in\cP_\sigma(\mu,\sigma)$,
\[
    \sum_{i=1}^N\bar w_i\bar p_i\bP(X>\bar p_i)
    =\sum_{j:\,w_j^*>0}\ \sum_{i\in C_j}\frac{w_j^*}{|C_j|}\,p_j\bP(X>p_j)
    =\sum_{j:\,w_j^*>0}w_j^*p_j\bP(X>p_j)
    =\sum_{j=1}^Nw_j^*p_j\bP(X>p_j),
\]
so that $R(\bar p)\ge R(p)$. Hence, $R(\bar p)=R(p)$, and $\bar w$ is optimal for $\bar p$.

Positive weights make the seller indifferent between all prices. By \eqref{eq:minimax}, $\max_ia^*_i=R(p)$, so that $a_i^*\le R(p)$ for all $i$, while $\sum_iw_i^*a_i^*=R(p)$. Therefore $a_i^*=R(p)$ for all $i$. As $a^*\in\bar{\mathcal A}(p)$, there are $\bP_n\in\cP_\sigma(\mu,\sigma)$ with $p_i\bP_n(X>p_i)\to R(p)$ for all $i$. With the intervals $I_1:=[0,p_1]$, $I_i:=(p_{i-1},p_i]$ for $2\le i\le N$ and $I_{N+1}:=(p_N,\infty)$, it follows that $\bP_n(I_i)\to m_i$, where
\begin{align}\label{eq:intervalmasses}
    m_1:=1-\frac{R(p)}{p_1},\qquad m_i:=R(p)\Bigl(\frac1{p_{i-1}}-\frac1{p_i}\Bigr)\ \ (2\le i\le N),\qquad m_{N+1}:=\frac{R(p)}{p_N}.
\end{align}
In particular, $m_{N+1}>0$, and a repeated price gives an empty interval with $m_i=0$.

We now locate the mass within each interval. Let $G(x):=\sum_{i=1}^Nw^*_ip_i\1_{\{x>p_i\}}$ denote the revenue of the menu $(w^*,p)$ at valuation $x$. We further let $\mathcal{M}(\R_+)$ be the set of non-negative measures defined on the measurable space $([0,\infty),\mathcal{B})$. The inner problem in \eqref{eq:Rp} is the primal
\begin{align}\label{eq:primal}
    R(p) = \inf_{\bP \in \mathcal{M}(\R_+)}\ \int G(x)\,d\bP(x)\quad\text{s.t.}\quad\int d\bP(x)=1,\quad\int x\,d\bP(x)=\mu,\quad\int x^2\,d\bP(x)=\mu^2+\sigma^2,
\end{align}
with corresponding dual
\begin{align}\label{eq:dual}
    \sup_{\lambda\in\R^3}\ \lambda_0+\lambda_1\mu+\lambda_2(\mu^2+\sigma^2)\qquad\text{s.t.}\qquad F(x):=\lambda_0+\lambda_1x+\lambda_2x^2\le G(x)\ \ \text{for all }x\ge0.
\end{align} The following lemma, proven in Appendix~\ref{ap:strong_duality}, shows that the primal and dual values in fact coincide.

\begin{lemma}[Strong duality]\label{lem:strong_dualityvar}
The optimal values of \eqref{eq:primal} and \eqref{eq:dual} coincide, and the supremum in \eqref{eq:dual} is attained.
\end{lemma}

Let $F^*(x)=\lambda_0^*+\lambda_1^*x+\lambda_2^*x^2$ attain the supremum in \eqref{eq:dual}. By Lemma~\ref{lem:strong_dualityvar}, for every $\bP\in\cP_\sigma(\mu,\sigma)$,
\[
    \int F^*(x)\,d\bP(x)=\lambda^*_0+\lambda_1^*\mu+\lambda_2^*(\mu^2+\sigma^2)=R(p).
\] We have $\lambda_2^*<0$. Otherwise, either $\lambda_2^*>0$, which implies that $F^*(x)$ is unbounded from above, contradicting $F^*(x)\le G(x)$ for sufficiently large $x$, or $\lambda_2^*=0$. In the latter case, if $\lambda^*_1 > 0$, we obtain the same contradiction, and therefore $\lambda_1^* \leq0$. But then $R(p)=\lambda_0^*+\lambda_1^*\mu\le\lambda_0^*=F^*(0)\le G(0)=0$, contradicting $R(p)>0$. Hence, $F^*$ is a strictly concave parabola, and we denote its vertex by $z(p)$.

We next show that $z(p)>p_N$. If $z(p)\le p_N$, then $F^*(x)$ is decreasing on $(p_N,\infty)$, so that for $x>p_N$, we have
\[
    G(x)-F^*(x)\geq G(x)-F^*(p_N)\geq G(x)-G(p_N)\geq w_N^*p_N>0.
\]
Since $G(x)-F^*(x)\ge0$ for all $x\ge0$, it follows that
\[
    \int(G(x)-F^*(x))\,d\bP_n(x)\ \ge\ \int_{I_{N+1}}(G(x)-F^*(x))\,d\bP_n(x)\ \ge\ w_N^*p_N\,\bP_n(I_{N+1})\ \longrightarrow\ w_N^*p_N\,m_{N+1}>0.
\]
However, the left-hand side converges to zero, because $\int G(x)\,d\bP_n(x)\to R(p)$ and $\int F^*(x)\,d\bP_n(x)=R(p)$. Hence, $z(p) > p_N$, and consequently, $F^*(x)$ is a strictly concave parabola with vertex $z(p)$, so that $$F^*(x)=F^*(z(p))+\lambda_2^*(x-z(p))^2.$$

Write $p_{N+1}:=z(p)$. We now locate Nature's mass within each interval using complementary slackness. Since $\int G(x)\,d\bP_n(x)\to R(p)$, while $\int F^*(x)\,d\bP_n(x)=R(p)$, we have
\begin{align}\label{eq:gap}
    \int(G(x)-F^*(x))\,d\bP_n(x)\longrightarrow0,
\end{align}
where $G(x)-F^*(x)\ge0$ for all $x\ge0$ by dual feasibility. In the limit, Nature's mass can therefore only remain at points $x$ where $G(x)=F^*(x)$. On $I_i$ with $i\le N$, $G(x)$ is constant while $F^*(x)$ is strictly increasing, so that $G(x)=F^*(x)$ is possible only at the right endpoint $x=p_i$. On $I_{N+1}$, $G(x)$ is constant as well, while $F^*(x)<F^*(z(p))$ for all $x\neq z(p)$, so that $G(x)=F^*(x)$ is possible only at the vertex $x=z(p)$.

A candidate for Nature's worst-case distribution is therefore
\begin{align}\label{eq:worstcase}
    \bP^*=\sum_{i=1}^{N+1}m_i\,\delta_{p_i},
\end{align}
which places the mass $m_i$ at $p_i$. For this candidate to be an actual worst-case distribution, it remains to show that $\bP^*$ satisfies the moment constraints
\begin{align}\label{eq:limitmoments}
    \sum_{i=1}^{N+1}m_ip_i=\mu\qquad\text{and}\qquad\sum_{i=1}^{N+1}m_ip_i^2=\mu^2+\sigma^2.
\end{align} To this end, we first show that
\begin{align}\label{eq:quadgap}
    G(x)-F^*(x)\ \ge\ -\lambda_2^*(x-p_i)^2\qquad\text{for all }x\in I_i,\quad i=1,\dots,N+1.
\end{align}
Recall that $F^*(x)=F^*(z(p))+\lambda_2^*(x-z(p))^2$. For $i\le N$ and $x\in I_i$, we have $G(x)=G(p_i)\ge F^*(p_i)$ and $x\le p_i\le z(p)$, so that
\begin{align*}
    G(x)-F^*(x)\ &\ge\ F^*(p_i)-F^*(x)=\lambda_2^*\bigl((p_i-z(p))^2-(x-z(p))^2\bigr)\\
    &=-\lambda_2^*(p_i-x)(2z(p)-p_i-x)\ \ge\ -\lambda_2^*(p_i-x)^2.
\end{align*} For $x\in I_{N+1}$, we have $G(x)=G(z(p))\ge F^*(z(p))$, because $z(p)>p_N$, so that $G(x)-F^*(x)\ge F^*(z(p))-F^*(x)=-\lambda_2^*(x-z(p))^2$. Hence, for $i=1,\dots,N+1$,
\[
    -\lambda_2^*\int_{I_i}(x-p_i)^2\,d\bP_n(x)\ \le\ \int(G(x)-F^*(x))\,d\bP_n(x)\ \longrightarrow\ 0.
\]
Since $\lambda_2^*<0$, it follows that
\[
    \int_{I_i}(x-p_i)^2\,d\bP_n(x)\to0.
\]

By the Cauchy--Schwarz inequality with $f(x)=x-p_i$ and $g(x)=1$,
\[
    \Bigl(\int_{I_i}(x-p_i)\,d\bP_n(x)\Bigr)^2\ \le\ \bP_n(I_i)\int_{I_i}(x-p_i)^2\,d\bP_n(x)\ \le\ \int_{I_i}(x-p_i)^2\,d\bP_n(x)\ \longrightarrow\ 0,
\]
so that also 
\[
    \int_{I_i}(x-p_i)\,d\bP_n(x)\to0.
\]

Therefore, since $\bP_n(I_i)\to m_i$, we obtain
\begin{align*}
    \int_{I_i}x\,d\bP_n(x)&=p_i\bP_n(I_i)+\int_{I_i}(x-p_i)\,d\bP_n(x)\ \to\ p_im_i,\\
    \int_{I_i}x^2\,d\bP_n(x)&=p_i^2\bP_n(I_i)+2p_i\int_{I_i}(x-p_i)\,d\bP_n(x)+\int_{I_i}(x-p_i)^2\,d\bP_n(x)\ \to\ p_i^2m_i.
\end{align*}

Summing over $i$ gives
\begin{align*}
    \mu &=\sum_{i=1}^{N+1}\int_{I_i}x\,d\bP_n(x)\ \longrightarrow\ \sum_{i=1}^{N+1}p_im_i,\\
    \mu^2+\sigma^2 &= \sum_{i=1}^{N+1}\int_{I_i}x^2\,d\bP_n(x)\ \longrightarrow\ \sum_{i=1}^{N+1}p_i^2m_i.
\end{align*}
Since $\mu$ and $\sigma$ do not depend on $n$, this proves \eqref{eq:limitmoments}, and consequently that $\bP^*\in\cP_\sigma(\mu,\sigma)$.

Next, we insert \eqref{eq:intervalmasses} into \eqref{eq:limitmoments}. Writing $r_i:=p_{i+1}/p_i$ for $i=1,\dots,N$, we obtain
\begin{align}\label{eq:uppermoments}
    \mu=p_1+R(p)\sum_{i=1}^N(r_i-1),\qquad\mu^2+\sigma^2=p_1^2+R(p)\sum_{i=1}^Np_i(r_i^2-1).
\end{align}

Multiplying the first equation by $2z(p)$, subtracting the second, and collecting the terms without $R(p)$ on the left-hand side gives
\[
    2z(p)(\mu-p_1)-(\mu^2+\sigma^2-p_1^2)=R(p)\sum_{i=1}^N\bigl(2z(p)(r_i-1)-p_i(r_i^2-1)\bigr)=R(p)\sum_{i=1}^N(r_i-1)(2z(p)-p_i-p_{i+1}),
\]
where the second equality uses $p_i(r_i^2-1)=(r_i-1)(p_i+p_{i+1})$, and the left-hand side equals $(z(p)-p_1)^2-(z(p)-\mu)^2-\sigma^2$. Hence,
\begin{align}\label{eq:Rformula}
    R(p) = \frac{(z(p)-p_1)^2-(z(p)-\mu)^2-\sigma^2}{\sum_{i=1}^N(r_i-1)(2z(p)-p_i-p_{i+1})}.
\end{align}

To determine the optimal prices, we would like to maximize $R(p)$ over $p$. However, this is not directly possible, as $R(p)$ depends on $z(p)$, which is known only implicitly. We therefore replace $z(p)$ in \eqref{eq:Rformula} by a free variable $z$ and maximize over $p$ and $z$ jointly. This gives an upper bound on $R^*$, which we will then show to be tight. For all $(p,z)$ with $0<p_1\le\dots\le p_N\le z$ and $p_1<z$, define
\begin{align}\label{eq:Psi}
    \Psi(p,z):=\frac{(z-p_1)^2-(z-\mu)^2-\sigma^2}{\sum_{i=1}^N(r_i-1)(2z-p_i-p_{i+1})},
\end{align}
where $p_{N+1}=z$ and $r_i=p_{i+1}/p_i$ for $i=1,\dots,N$. Additionally, denote the numerator and the denominator in \eqref{eq:Psi} by $A(p,z)$ and $D(p,z)$, respectively. Notice that $D(p,z)>0$, because $p_1\le\dots\le p_N\le z$ causes every term of $D(p,z)$ to be non-negative, and since $p_1<z$, at least one term is positive. With $\Psi^*:=\sup_{p,z}\Psi(p,z)$, we have
\begin{align}\label{eq:upper}
    0<R^*\le\Psi^*.
\end{align}

To determine $\Psi^*$, we proceed in three steps. First, we show that the supremum is attained. Second, we show that every maximizer is an interior point of the domain, so that it satisfies the first-order conditions. Third, we show that the first-order conditions have at most one solution, which therefore determine $\Psi^*$. We start with the first step.

It suffices to show that for every $c\in(0,\Psi^*)$, the set of all $(p,z)$ with $\Psi(p,z)\geq c$ is compact. For this, we show that $p_1$ is bounded away from zero, $z$ is bounded away from $p_1$, and $z$ is bounded from above. Let $\Psi(p,z)\geq c$, so that $A(p,z)\geq cD(p,z)>0$. Then $(z-p_1)^2>(z-\mu)^2+\sigma^2\ge\sigma^2$, so that $z-p_1>\sigma$. Moreover, $p_1<\mu$, as otherwise, $\mu-p_1\le0$ and $2z-\mu-p_1>0$, because $z>p_1\geq\mu$, so that $A(p,z)=(\mu-p_1)(2z-\mu-p_1)-\sigma^2<0$, which contradicts $A(p,z)>0$. As a result, $\mu^2+\sigma^2-p_1^2>0$, so that $A(p,z)=2z(\mu-p_1)-(\mu^2+\sigma^2-p_1^2)<2z(\mu-p_1)<2\mu z$. Together with $A(p,z)\ge cD(p,z)$, this gives $D(p,z)<2\mu z/c$.

Next, we show that
\[
    D(p,z)\ \geq\ \frac z2\ln\Bigl(\frac z{2p_1}\Bigr).
\]
If $p_1>z/2$, this holds because the right-hand side is negative. Otherwise, let $k$ be the largest index with $p_k\leq z/2$, so that $p_{k+1}>z/2$. For $i\leq k$, we have $2z-p_i-p_{i+1}\geq z/2$ and $r_i-1\geq \ln r_i$, since $\ln x \leq x - 1$ for all $x >0$, and therefore
\[
    D(p,z)\ \geq\ \frac z2\sum_{i=1}^k\ln r_i=\frac z2\ln\Bigl(\frac{p_{k+1}}{p_1}\Bigr)\ \geq\ \frac z2\ln\Bigl(\frac z{2p_1}\Bigr).
\]
Combining the two bounds on $D(p,z)$ gives $z<2e^{4\mu/c}p_1$. With $p_1<\mu$ and $z>z-p_1>\sigma$, we obtain $z<2\mu e^{4\mu/c}$ and $p_1>\sigma e^{-4\mu/c}/2$. Together with $z-p_1>\sigma$, this shows that the set is compact, and consequently that $\Psi^*$ is attained at some $(p^*,z^*)$.

Next, we show that every maximizer is an interior point of the domain, i.e., has distinct entries $0<p_1<\dots<p_N<z$. Suppose that a maximizer $(p,z)$ has a repeated entry $p_{k-1}=p_k$ for some $2\le k\le N+1$. Recall that, with $p_{N+1}=z$,
\[
    D(p,z)=\sum_{i=1}^N(r_i-1)(2z-p_i-p_{i+1})=\sum_{i=1}^N\frac{(p_{i+1}-p_i)(2z-p_i-p_{i+1})}{p_i}=\sum_{i=1}^N\frac{(z-p_i)^2-(z-p_{i+1})^2}{p_i},
\]
so that each term depends only on the pair $(p_{i},p_{i+1})$, and is zero if $p_i=p_{i+1}$. Let $p'$ be obtained by removing $p_{k-1}$ and adding a copy of $p_1$, so that $p'=(p_1,p_1,p_2,\dots,p_{k-2},p_k,\dots,p_N)$. Then $D(p,z)=D(p',z)$. Indeed, both sums contain one zero term, for the pair $(p_{k-1},p_k)$ and the pair $(p'_1,p'_2)=(p_1,p_1)$, respectively, while all remaining terms coincide. Moreover, $A(p,z)$ depends only on $p_1$ and $z$, so that $A(p,z)=A(p',z)$ and $\Psi(p,z)=\Psi(p',z)$. Hence, $(p',z)$ is also a maximizer.

Repeating this, we may move all repeated entries to copies of $p_1$, and we denote the resulting vector again by $p'$. Then $(p',z)$ is a maximizer with $p'_{i-1}=p'_i<p'_{i+1}$ for some $2\leq i\leq N$, where $p'_i$ is the last copy of $p_1$. Next, observe that
\begin{align}\label{eq:partials}
    \frac{\partial\Psi(p,z)}{\partial p_i}
    &=-\frac{A(p,z)}{D(p,z)^2}\,\frac{\partial D(p,z)}{\partial p_i}\nonumber\\
    &=-\frac{A(p,z)}{p_i^2D(p,z)^2}\biggl[\frac{2p_i^2(z-p_i)}{p_{i-1}}-2p_i(z-p_i)-(z-p_i)^2+(z-p_{i+1})^2\biggr]\nonumber\\
    &=-\frac{A(p,z)}{p_i^2D(p,z)^2}\Bigl[(z-p_{i+1})^2-(z-p_i)(z+p_i-2r_{i-1}p_i)\Bigr].
\end{align}
Since $A(p',z)>0$, it holds that $\partial\Psi(p',z)/\partial p_i > 0$ if and only if $$(z-p'_{i+1})^2-(z-p'_i)(z+p'_i-2(p'_i/p'_{i-1})p'_i) = (z - p'_{i+1})^2 - (z - p'_{i})^2 < 0,$$ where the equality uses $p'_{i-1}=p'_i$. Furthermore, the inequality holds because $p'_i<p'_{i+1}\leq z$. As $p'_i<p'_{i+1}$, increasing $p'_i$ slightly keeps $(p',z)$ in the domain and increases $\Psi(p',z)$, which contradicts maximality. Therefore, every maximizer satisfies $0<p^*_1<\dots<p^*_N<z^*$, i.e., is an interior point of the domain.

Let $(p^*,z^*)$ be any maximizer. We now evaluate the first-order conditions. Since $A(p^*,z^*)>0$, \eqref{eq:partials} shows that $\partial\Psi(p^*,z^*)/\partial p_i=0$ if and only if
\begin{align}\label{eq:sqrd_eq}
    (z^*-p^*_{i+1})^2=(z^*-p^*_i)(z^*+p^*_i-2r^*_{i-1}p^*_i),\qquad i=2,\dots,N.
\end{align}
For $2 \leq i\leq N-1$, taking square roots and using $p^*_{i+1}<z^*$ gives \eqref{eq:rec}, and $i=N$ gives \eqref{eq:terminal}, because the left-hand side is zero and $z^*>p^*_N$.

For the derivative with respect to $z$, note that $\partial A(p,z)/\partial z=2(\mu-p_1)$ and $\partial D(p,z)/\partial z=2\sum_{i=1}^N(r_i-1)$, so that
\[
    \frac{\partial\Psi(p,z)}{\partial z}
    =\frac{2(\mu-p_1)D(p,z)-2A(p,z)\sum_{i=1}^N(r_i-1)}{D(p,z)^2}
    =\frac{2\sum_{i=1}^N(r_i-1)}{D(p,z)}\biggl[\frac{\mu-p_1}{\sum_{i=1}^N(r_i-1)}-\Psi(p,z)\biggr].
\]
Since $\sum_{i=1}^N(r^*_i-1)>0$ and $D(p^*,z^*)>0$, as $p^*_1<\dots<p^*_N<z^*$, it follows that $\partial\Psi(p^*,z^*)/\partial z=0$ if and only if
\begin{align}\label{eq:foc_z}
    \Psi^*=\frac{\mu-p^*_1}{\sum_{i=1}^N(r^*_i-1)}.
\end{align}

Using $(r_i-1)(p_i+p_{i+1})=p_i(r_i^2-1)$, we may write
\begin{align}\label{eq:foc_z2}
    \Psi^*=\frac{A(p^*,z^*)}{D(p^*,z^*)}=\frac{2z^*(\mu-p^*_1)-\bigl(\mu^2+\sigma^2-(p^*_1)^2\bigr)}{2z^*\sum_{i=1}^N(r^*_i-1)-\sum_{i=1}^Np^*_i\bigl((r^*_i)^2-1\bigr)}.
\end{align}
Equating \eqref{eq:foc_z2} with \eqref{eq:foc_z} gives \eqref{eq:moment}.

Finally, consider the derivative with respect to $p_1$. Writing
\begin{align*}
    D(p,z)&=\sum_{i=1}^N\frac{(z-p_i)^2-(z-p_{i+1})^2}{p_i}=\sum_{i=1}^N\frac{(2zp_{i+1}-p_{i+1}^2)-(2zp_i-p_i^2)}{p_i}\\
    &=\sum_{i=1}^N\biggl(\frac{2zp_{i+1}-p_{i+1}^2}{p_i}-2z+p_i\biggr),
\end{align*}
we obtain
\[
    \frac{\partial D(p,z)}{\partial p_1}=1-\frac{2zp_2-p_2^2}{p_1^2},\qquad
    \frac{\partial D(p,z)}{\partial p_j}=\frac{2(z-p_j)}{p_{j-1}}+1-\frac{2zp_{j+1}-p_{j+1}^2}{p_j^2},\quad j=2,\dots,N.
\]
By \eqref{eq:partials} and \eqref{eq:sqrd_eq}, the derivatives with respect to $p_2,\dots,p_N$ are zero at $(p^*,z^*)$. Multiplying them by $p^*_j$ and summing over $j=2,\dots,N$ gives
\begin{align*}
    0&=\sum_{j=2}^N\biggl(\frac{2z^*p^*_j-(p^*_j)^2}{p^*_{j-1}}-\frac{2z^*p^*_{j+1}-(p^*_{j+1})^2}{p^*_j}\biggr)-\sum_{j=2}^N\biggl(\frac{(p^*_j)^2}{p^*_{j-1}}-p^*_j\biggr)\\
    &=\frac{2z^*p^*_2-(p^*_2)^2}{p^*_1}-\frac{(z^*)^2}{p^*_N}-\sum_{j=2}^N\biggl(\frac{(p^*_j)^2}{p^*_{j-1}}-p^*_j\biggr)\\
    &=\frac{2z^*p^*_2-(p^*_2)^2}{p^*_1}-p^*_1-\sum_{i=1}^Np^*_i\bigl((r^*_i)^2-1\bigr),
\end{align*}
so that
\[
    \frac{\partial D(p^*,z^*)}{\partial p_1}=-\frac{1}{p^*_1}\sum_{i=1}^Np^*_i\bigl((r^*_i)^2-1\bigr).
\]
Moreover, combining \eqref{eq:foc_z} with \eqref{eq:moment} gives
\[
    \Psi^*\sum_{i=1}^Np^*_i\bigl((r^*_i)^2-1\bigr)=\mu^2+\sigma^2-(p^*_1)^2.
\]
Since $\partial A(p,z)/\partial p_1=-2(z-p_1)$, the quotient rule therefore gives
\[
    \frac{\partial\Psi(p^*,z^*)}{\partial p_1}=\frac{1}{D(p^*,z^*)}\biggl[\frac{\partial A(p^*,z^*)}{\partial p_1}-\Psi^*\,\frac{\partial D(p^*,z^*)}{\partial p_1}\biggr]=\frac{\mu^2+\sigma^2-(p^*_1)^2-2p^*_1(z^*-p^*_1)}{p^*_1D(p^*,z^*)},
\]
which is zero if and only if \eqref{eq:zdef} holds.

Every maximizer is therefore generated by its first two prices through \eqref{eq:zdef} and \eqref{eq:rec}, and it satisfies $0<p^*_1<\dots<p^*_N$ together with \eqref{eq:terminal} and \eqref{eq:moment}. By the following lemma, at most one pair of prices has this property, so that the maximizer $(p^*,z^*)$ is unique. The proof is found in Appendix~\ref{ap:Lemma_var}.

\begin{lemma}[Uniqueness]\label{lem:e_and_u}
Let $N\ge2$. At most one pair $(p_1,p_2)$ generates, through \eqref{eq:zdef} and \eqref{eq:rec}, a vector $0<p_1<\dots<p_N$ that satisfies \eqref{eq:terminal} and \eqref{eq:moment}.
\end{lemma}

It remains to show that $\Psi^*\leq R^*$ and to determine $w^*$. If we knew that $z^*=z(p^*)$, tightness would follow from \eqref{eq:Rformula}, since then $\Psi^*=\Psi(p^*,z(p^*))=R(p^*)\leq R^*$. However, $z(p^*)$ is known only implicitly, so that we cannot verify this directly. Suppose, however, that the upper bound \eqref{eq:upper} \textit{were} tight, so that $z(p^*)=z^*$. The analysis for fixed prices would then determine $w^*$ for the prices $p^*$, and we can use these weights to construct a candidate menu. We then verify, without assuming tightness, that this menu guarantees $\Psi^*$.

Indeed, the dual parabola $F^*(x)=F^*(z^*)+\lambda_2^*(x-z^*)^2$ would then coincide with $G(x)$ at $p^*_1,\dots,p^*_N$ and $z^*$. Since $G(x)$ increases by $w^*_ip^*_i$ between $p^*_i$ and $p^*_{i+1}$, with $p^*_{N+1}=z^*$, it follows that $w^*_ip^*_i=F^*(p^*_{i+1})-F^*(p^*_i)=-\lambda_2^*\bigl((z^*-p^*_i)^2-(z^*-p^*_{i+1})^2\bigr)$, and since the weights have to sum to one, $\lambda_2^*=-1/D(p^*,z^*)$. The candidate weights are therefore
\begin{align}\label{eq:weights_var}
    w_i=\frac{(z^*-p^*_i)^2-(z^*-p^*_{i+1})^2}{p^*_iD(p^*,z^*)}=\frac{(r^*_i-1)(2z^*-p^*_i-p^*_{i+1})}{\sum_{j=1}^N(r^*_j-1)(2z^*-p^*_j-p^*_{j+1})},\qquad i=1,\dots,N,
\end{align}
which are the weights \eqref{eq:var_menu}.

To verify that the menu $(w,p^*)$ guarantees $\Psi^*$, consider
\begin{align}\label{eq:primal_w}
    \inf_{\bP \in \mathcal{M}(\R_+)}\ \int G_w(x)\,d\bP(x)\quad\text{s.t.}\quad\int d\bP(x)=1,\quad\int x\,d\bP(x)=\mu,\quad\int x^2\,d\bP(x)=\mu^2+\sigma^2,
\end{align}
where $G_w(x):=\sum_{i=1}^Nw_ip^*_i\1_{\{x>p^*_i\}}$, and corresponding dual
\begin{align}\label{eq:dual_w}
    \sup_{\lambda\in\R^3}\ \lambda_0+\lambda_1\mu+\lambda_2(\mu^2+\sigma^2)\qquad\text{s.t.}\qquad F_w(x):=\lambda_0+\lambda_1x+\lambda_2x^2\le G_w(x)\ \ \text{for all }x\ge0.
\end{align}

Consider the dual solution
\[
    \lambda_0=\frac{(p^*_1)^2-2z^*p^*_1}{D(p^*,z^*)},\qquad\lambda_1=\frac{2z^*}{D(p^*,z^*)},\qquad\lambda_2=-\frac{1}{D(p^*,z^*)},
\]
for which
\[
    F_w(x)=\lambda_0+\lambda_1x+\lambda_2x^2=\frac{(z^*-p^*_1)^2-(x-z^*)^2}{D(p^*,z^*)}.
\]
Its objective value equals
\[
    \lambda_0+\lambda_1\mu+\lambda_2(\mu^2+\sigma^2)=\frac{2z^*(\mu-p^*_1)-\bigl(\mu^2+\sigma^2-(p^*_1)^2\bigr)}{D(p^*,z^*)}=\frac{A(p^*,z^*)}{D(p^*,z^*)}=\Psi^*.
\]
If this solution is feasible for \eqref{eq:dual_w}, weak duality gives
\[
    \Psi^*\ \le\ \inf_{\bP\in\cP_\sigma(\mu,\sigma)}\int G_w(x)\,d\bP(x)\ \le\ R(p^*)\ \le\ R^*.
\]
It therefore remains to show that $F_w(x)\le G_w(x)$ for all $x\ge0$.

First, $F_w(x)=G_w(x)$ for $x\in\{p^*_1,\dots,p^*_N,z^*\}$. Indeed, with $p^*_{N+1}=z^*$, we have $F_w(p^*_1)=0=G_w(p^*_1)$, and by \eqref{eq:weights_var},
\[
    F_w(p^*_{i+1})-F_w(p^*_i)=\frac{(z^*-p^*_i)^2-(z^*-p^*_{i+1})^2}{D(p^*,z^*)}=w_ip^*_i=G_w(p^*_{i+1})-G_w(p^*_i),\qquad i=1,\dots,N.
\]
Summing these differences gives, for $i=1,\dots,N+1$,
\[
    F_w(p^*_i)=F_w(p^*_1)+\sum_{j=1}^{i-1}\bigl(F_w(p^*_{j+1})-F_w(p^*_j)\bigr)=G_w(p^*_1)+\sum_{j=1}^{i-1}\bigl(G_w(p^*_{j+1})-G_w(p^*_j)\bigr)=G_w(p^*_i).
\]
Second, $F_w(x)$ is increasing on $[0,z^*]$ and maximal at $z^*$, while $G_w(x)$ is constant on each of the intervals $[0,p^*_1]$, $(p^*_{i-1},p^*_i]$ for $2\le i\le N$, and $(p^*_N,\infty)$. On each of these intervals, $F_w(x)$ is therefore bounded by its value at the right endpoint, or at $z^*$ for the last interval, where it equals $G_w(x)$. Hence, $F_w(x)\le G_w(x)$ for all $x\ge0$.

Hence, $R^*\geq\Psi^*$. Together with \eqref{eq:upper}, it follows that $R^*=\Psi^*$, so that the menu $(w^*,p^*)$ of Theorem~\ref{th:var_N} is optimal. This finishes the proof.

\section{Proof of the Mean--MAD Finite Menu Theorem}\label{sec:bijection}
The proof of Theorem~\ref{th:mad_N} rests on a structural result that is of independent interest and that explains why mean--MAD ambiguity is the more tractable of the two: for a broad class of objectives, the worst case under mean--MAD ambiguity coincides with the worst case under mean--support ambiguity, in which the dispersion constraint is replaced by an upper bound on the valuation. We state this result first, and return to the proof of Theorem~\ref{th:mad_N} at the end of the section.
 
\begin{theorem}\label{th:bijection}
    Let $G: [0,\infty) \to \R$ be non-decreasing and differentiable for all sufficiently large $x$, with $\lim_{x\to\infty}G'(x)=0$. For $\mu>0$ and $d\in(0,2\mu)$,
    \begin{align}\label{eq:bijection_mad}
    \inf_{\bP\in\cP_d(\mu,d)}\bE_{\bP}[G(X)] = \inf_{\bP\in\cP_{\beta}(\mu-d/2,\mu)}\bE_{\bP}[G(X)].
    \end{align}
    Equivalently, for $\beta>0$ and $m\in(0,\beta)$,
    \begin{align}\label{eq:bijection_range}
    \inf_{\bP\in\cP_{\beta}(m,\beta)}\bE_{\bP}[G(X)] = \inf_{\bP\in\cP_{d}(\beta,2(\beta-m))}\bE_{\bP}[G(X)].
    \end{align}
\end{theorem}
 
The proof of Theorem~\ref{th:bijection} is given in Section~\ref{sec:proof}. To see where this correspondence comes from, we begin with an elementary identity. Let $X$ be any integrable random variable with mean $\mu$. Then
\begin{align}\label{eq:trunc_identity}
    \bE[\min\{X,\mu\}] = \mu - \tfrac{1}{2}\bE|X-\mu|.
\end{align}
Indeed, $\min\{X,\mu\} = \mu - (\mu - X)^+$, where $(\mu - X)^+ = \max\{\mu - X,0\}$. Taking expectations yields $\bE[\min\{X,\mu\}] = \mu - \bE[(\mu - X)^+]$. Moreover, since $|X-\mu| = (X-\mu)^+ + (\mu - X)^+$ and $\bE[X-\mu]=0$, we have $\bE[(X-\mu)^+] = \bE[(\mu - X)^+]$, implying $\bE|X-\mu| = 2\bE[(\mu - X)^+]$. Substituting gives \eqref{eq:trunc_identity}.
This identity shows that the mean absolute deviation can be interpreted as a two-sided dispersion around $\mu$, which, upon truncation at $\mu$, translates into a downward shift in the mean. The truncation caps values at $\mu$, so a distribution with mean--MAD parameters $(\mu,d)$ is mapped to one supported on $[0,\mu]$, while the identity fixes its mean at $\mu-d/2$. This explains the parameters in \eqref{eq:bijection_mad}, but not the equality. That the worst cases nevertheless coincide is the content of Theorem~\ref{th:bijection}, and this is where the conditions on $G$ enter.
 
These conditions are satisfied by the expected revenue $\sum_{i=1}^Nw_ip_i\1\{x\geq p_i\}$ that a menu $(w,p)$ collects from a buyer with valuation $x$, which is non-decreasing and constant above $p_N$. The mean--MAD pricing problem therefore reduces to its mean--support counterpart, which is solved by \citet{wang2025power}:

\begin{proposition}[\citealp{wang2025power}, Proposition 11]\label{prop:wang}
For $\cP = \cP_\beta(m,\beta)$ and $N\in\N$, an optimal menu $(w^*,p^*)$ is given by $w^*_i = 1/N$ and $p^*_i = p^*_1(\beta/p^*_1)^{(i-1)/N}$, $i=1,\dots,N$, where $p^*_1$ is the unique solution to $m/p^*_1 + N(p^*_1/\beta)^{1/N} = N+1$. The optimal expected revenue is $(m-p^*_1)\,N^{-1}\big((\beta/p^*_1)^{1/N}-1\big)^{-1}$.
\end{proposition}
 
We now turn to the proof of Theorem \ref{th:mad_N}. Fix a menu $(w,p)\in\mathcal H$. Its expected revenue under $\bP$ equals $\bE_{\bP}[G(X)]$ with $G(x) = \sum_{i=1}^N w_ip_i\1\{x\geq p_i\}$, which is non-decreasing and constant for $x\geq p_N$. By \eqref{eq:bijection_mad}, the worst-case revenue of the menu is therefore the same under $\cP_d(\mu,d)$ and under $\cP_\beta(\mu-d/2,\mu)$. Since this holds for every menu,
\begin{align}\label{eq:bijection_result}
    R^* = \sup_{(w,p)\in\mathcal H}\inf_{\bP\in\cP_\beta(\mu-d/2,\mu)} \sum_{i=1}^Nw_ip_i\bP(X\geq p_i),
\end{align}
and the two problems have the same optimal menus. The right-hand side of \eqref{eq:bijection_result} is the finite-menu pricing problem under mean--support ambiguity with mean $\mu - d/2$ and support upper bound $\mu$. Applying Proposition~\ref{prop:wang} with $m = \mu - d/2$ and $\beta = \mu$ yields the menu \eqref{eq:mad_menu} and the equation \eqref{eq:p1}. Substituting \eqref{eq:p1} into the revenue formula of Proposition~\ref{prop:wang} gives \eqref{eq:rev_N_2a}. This finishes the proof.

\section{Proof of the MAD-Support Equivalence Theorem}\label{sec:proof}
We now turn to the proof of Theorem \ref{th:bijection}. Consider $\cP = \cP_d(\mu,d)$ and recall that $\mathcal{M}(\R_+)$ denotes the set of non-negative measures defined on the measurable space $([0,\infty),\mathcal{B})$. We can represent \eqref{eq:bijection_mad} as the following linear program:
\begin{equation}
\begin{aligned}\label{eq:LP1}
\! v_p = &\inf_{\mathbb{P}\in \mathcal{M}(\R_+)} &  &\int_x G(x){\rm d} \mathbb{P}(x),\\
&\text{s.t.} &      &  \int_x {\rm d}\mathbb{P}(x)=1, \ \int_x x{\rm d}\mathbb{P}(x)=\mu,\ \int_x |x-\mu|{\rm d}\mathbb{P}(x) = d.
\end{aligned}
\end{equation}
Consider the corresponding dual problem
\begin{equation}\label{eq:LP2}
\begin{aligned}
v_d = &\sup_{\lambda \in \mathbb{R}^3}
    \big[\,\lambda_0 + \lambda_1 \mu+\lambda_2 d\,\big],\\
&\text{s.t.}\quad 
    G(x) 
    \;\geq\;
    \lambda_0 + \lambda_1 x+\lambda_2 |x-\mu| 
    \;=:\; F(x;\lambda), \quad \forall x\in[0,\infty).
\end{aligned}
\end{equation} 

We will first establish that $v_p = v_d$.
\begin{lemma}\label{lem:strong_duality}
Strong duality holds for \eqref{eq:LP1} and \eqref{eq:LP2}, i.e., $v_p = v_d$.
\end{lemma}
The proof is given in Appendix \ref{ap:strong_duality}. Next, we will argue that no feasible solution exists that satisfies $\lambda_1+\lambda_2>0$. Assume $\lambda = (\lambda_0,\lambda_1,\lambda_2)$ is feasible with $\lambda_1+\lambda_2 > 0$. Since $\lim_{x\xrightarrow{}\infty}G'(x)=0$, for every $\epsilon>0$ there exists $X_{\epsilon}$ such that $G$ is differentiable on $(X_\epsilon,\infty)$ and $G'(x)\leq\epsilon$ for all $x\geq X_\epsilon$. Hence, for every $x>X_\epsilon$, the mean value theorem gives some $\xi\in(X_\epsilon,x)$ such that
\begin{align*}
    G(x)-G(X_\epsilon)
    =G'(\xi)(x-X_\epsilon)
    \leq\epsilon(x-X_\epsilon).
\end{align*}
Take $\epsilon=(\lambda_1+\lambda_2)/2 > 0$ and $x > \max\{X_\epsilon,\mu\}$. Now observe that
\begin{align*}
    F(x;\lambda) = \lambda_0 + \lambda_1x + \lambda_2(x-\mu)
    = \lambda_0-\lambda_2\mu + (\lambda_1+\lambda_2)x
    \leq G(x) \leq G(X_\epsilon) +\epsilon(x-X_\epsilon),
\end{align*}
which implies that feasibility requires
\begin{align}\label{eq:>inf}
    \frac{\lambda_1+\lambda_2}{2}x \leq G(X_\epsilon) - [\lambda_0-\lambda_2\mu] - \frac{\lambda_1+\lambda_2}{2}X_\epsilon
\end{align}
to hold for all $x > \max\{X_\epsilon,\mu\}$. However, since the right-hand side of \eqref{eq:>inf} does not depend on $x$ and the left-hand side is a linear function with a positive slope, any $x$ sufficiently large will violate \eqref{eq:>inf}. This is a contradiction. Hence, any feasible solution satisfies $\lambda_1+\lambda_2\leq 0$.

We will now argue that one can restrict to feasible solutions that satisfy $\lambda_1+\lambda_2=0$. Assume $\lambda=(\lambda_0,\lambda_1,\lambda_2) \in \R^3$ is feasible and satisfies $\lambda_1+\lambda_2 < 0$, then there exists a $\delta>0$ such that $\lambda_1+\lambda_2 + \delta = 0$. Now take $\lambda'=(\lambda_0-\delta\mu,\lambda_1+\delta,\lambda_2)$ and notice that $\lambda'_1 + \lambda'_2 = 0$. Furthermore, $\lambda'$ is dual feasible, as for any $x \leq \mu$, it holds that $$F(x;\lambda')=\lambda_0-\delta\mu+(\lambda_1+\delta)x+\lambda_2(\mu-x) = F(x;\lambda)-\delta(\mu-x) \leq F(x;\lambda) \leq G(x),$$ and for any $x > \mu$ it holds that
\begin{align*}
    F(x;\lambda') &= \lambda_0-\delta\mu+(\lambda_1+\delta)x+\lambda_2(x-\mu) \\
    &= \lambda_0-\delta\mu+(\lambda_1+\delta)x-(\lambda_1+\delta)(x-\mu)\\ &= \lambda_0+\lambda_1\mu = F(\mu;\lambda) \leq G(\mu) \leq G(x).
\end{align*}
Furthermore, $\lambda_0-\delta\mu+(\lambda_1+\delta)\mu+\lambda_2d=\lambda_0+\lambda_1\mu+\lambda_2d$. Hence, $\lambda'$ attains the same dual objective value as $\lambda$. Therefore, it is without loss of objective value to restrict attention to feasible solutions satisfying $\lambda_1+\lambda_2=0$. 

For any such feasible solution, define $\bar{\lambda}_0 = \lambda_0-\lambda_1\mu$ and $\bar{\lambda}_1 = 2\lambda_1$. Then, every feasible solution of \eqref{eq:LP2} corresponds to a feasible solution of
\begin{equation}\label{eq:dual_range2}
\begin{aligned}
\bar{v}_d = &\sup_{\bar{\lambda} \in \mathbb{R}^2}
    \big[\,\bar{\lambda}_0+\bar{\lambda}_1(\mu-d/2)\,\big],\\
&\text{s.t.}\quad 
    G(x) 
    \;\geq\;
    \bar{\lambda}_0  + \bar{\lambda}_1 x
    \;=\; \bar{F}(x;\bar{\lambda}), \quad \forall x\in[0,\mu].
\end{aligned}
\end{equation} 

Consequently, $v_d \leq \bar{v}_d$. Conversely, take any feasible solution $\bar{\lambda}=(\bar{\lambda}_0,\bar{\lambda}_1)$ of \eqref{eq:dual_range2}. Define $\lambda_0 = \bar{\lambda}_0+(\bar{\lambda}_1\mu)/2$, $\lambda_1 = \bar{\lambda}_1/2$, and $\lambda_2 = -\bar{\lambda}_1/2$. We show that this gives a feasible solution of \eqref{eq:LP2}. If $x \leq \mu$, then
\begin{align*}
    \lambda_0+\lambda_1x+\lambda_2(\mu-x) = \bar{\lambda}_0+\bar{\lambda}_1x = \bar{F}(x;\bar{\lambda}) \leq G(x).
\end{align*}
Similarly, if $x > \mu$, then
\begin{align*}
    \lambda_0+\lambda_1x+\lambda_2(x-\mu) = \bar{\lambda}_0+ \bar{\lambda}_1\mu = \bar{F}(\mu;\bar{\lambda}) \leq G(\mu) \leq G(x).
\end{align*}
Hence, the constructed $\lambda$ is feasible for \eqref{eq:LP2}, and $\lambda_0+\lambda_1\mu+\lambda_2d=\bar{\lambda}_0+\bar{\lambda}_1\left(\mu-d/2\right)$, which is the objective value of \eqref{eq:dual_range2}. Therefore, $v_d \geq \bar{v}_d$, and hence, $v_d = \bar{v}_d$.

Next, define $\bar{\mathcal{M}}([0,\mu])$ as the set of all non-negative measures defined on the measurable space $([0,\mu],\mathcal{B})$. Then the corresponding primal problem of \eqref{eq:dual_range2} can be written as
\begin{equation}
\begin{aligned}\label{eq:primal_range}
\! \bar{v}_p = &\inf_{\mathbb{P}\in \bar{\mathcal{M}}([0,\mu])} &  &\int_x G(x){\rm d} \mathbb{P}(x),\\
&\text{s.t.} &      &  \int_x {\rm d}\mathbb{P}(x)=1, \ \int_x x{\rm d}\mathbb{P}(x)=\mu-\frac{d}{2}.
\end{aligned}
\end{equation}

We now establish that $\bar{v}_p = \bar{v}_d$.
\begin{lemma}\label{lem:strong_duality2}
Strong duality holds for \eqref{eq:dual_range2} and \eqref{eq:primal_range}, i.e., $\bar{v}_p = \bar{v}_d$.
\end{lemma}
The proof is given in Appendix \ref{ap:strong_duality}. Consequently, $v_p = v_d = \bar{v}_d = \bar{v}_p$, which proves \eqref{eq:bijection_mad}. The identity \eqref{eq:bijection_range} follows by applying \eqref{eq:bijection_mad} with $\mu = \beta$ and $d = 2(\beta-m)$. This finishes the proof.

The algebraic proof above establishes the equivalence. To complement it, we close the section with a geometric reading of the argument, which can be understood through three steps. \begin{figure}[h!]
  \centering
  \begin{subfigure}[b]{0.32\linewidth}
    \centering
    \caption*{\fontsize{9.5}{11}\selectfont Step 1: Any $F$ is affine then flat.}
    \includegraphics[width=\linewidth]{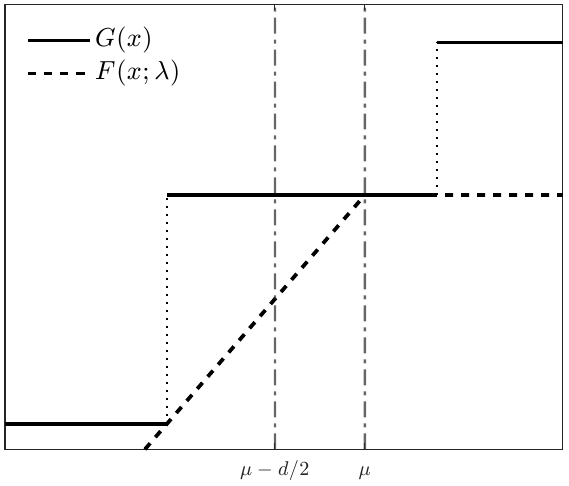}
  \end{subfigure}\hfill
  \begin{subfigure}[b]{0.32\linewidth}
    \centering
    \caption*{\fontsize{9.5}{11}\selectfont Step 2: Flat part adds no constraint.}
    \includegraphics[width=\linewidth]{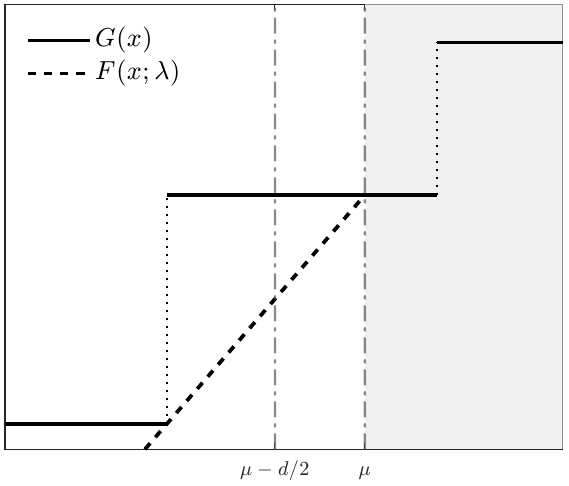}
  \end{subfigure}\hfill
  \begin{subfigure}[b]{0.32\linewidth}
    \centering
    \caption*{\fontsize{9.5}{11}\selectfont Step 3: This is mean--support.}
    \includegraphics[width=\linewidth]{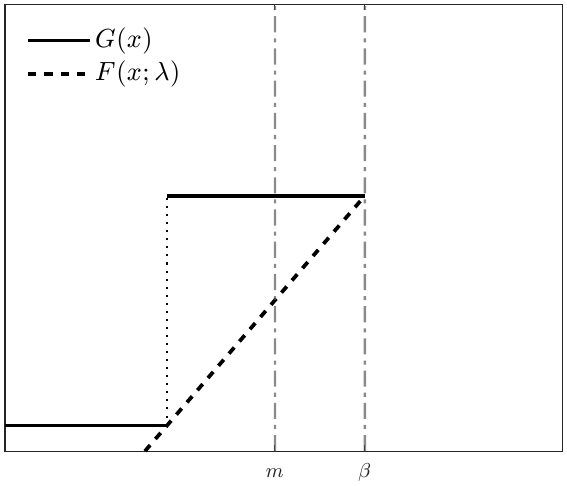}
  \end{subfigure}

  \caption{From mean--MAD to mean--support in three steps.}
  \label{fig:minorants}
\end{figure} First, in the mean--MAD setting the dual minorant $F(x;\lambda)$ is piecewise linear with a kink at $x=\mu$, induced by the absolute-value structure of the MAD constraint. Moreover, we can restrict $F$ to be non-decreasing affine on $[0,\mu]$ and constant on $[\mu,\infty)$ when $G$ satisfies the conditions of Theorem \ref{th:bijection}. Second, the constant part imposes no constraint. Since $G$ is non-decreasing, the constant part of $F$ is below $G$ for all $x>\mu$ when it is feasible at $x=\mu$. Third, what remains is a non-decreasing affine function on $[0,\mu]$, which is structurally identical to the dual minorant of a mean--support problem with support upper bound $\beta=\mu$ and mean $m = \mu - d/2$. Figure \ref{fig:minorants} illustrates these three steps. Conversely, starting from the mean--support setting, analogous steps construct a mean--MAD setting with $\mu = \beta$ and $d = 2(\beta-m)$, by extending the minorant with a redundant constant part on $[\beta,\infty)$.

\section{Conclusion}\label{sec:conclusion}
A seller with limited knowledge of demand must decide which summary statistics to rely on, and our results show that the choice of dispersion measure can have important consequences. As long as the seller posts a single price, the variance and the MAD give similar advice, which is to price lower in more volatile markets. However, once the seller offers a menu of more than one price, the two measures lead to qualitatively different outcomes. A seller who relies on the variance eventually raises the expected offered price, whereas a seller who relies on the MAD continues to lower all prices. The dispersion measure is therefore more than a matter of modeling convenience, as it carries an assumption about how the high valuations in the market are distributed, and this assumption remains hidden until the selling mechanism is sufficiently rich.

Our results also point to several directions for future research. One direction is to study finite menus under alternative objectives. In particular, the corresponding problems for the competitive ratio and absolute regret remain open under both mean--variance and mean--MAD ambiguity. Another direction is to apply the equivalence between mean--MAD and mean--support ambiguity outside of pricing. The newsvendor problem is a natural example of this, as sales are capped by the order quantity in the same way that revenue is capped by the price. In any such problem, a result obtained under one of the two ambiguity sets carries over to the other, in either direction. We expect the same to hold in other settings where a decision limits how much of a high realization can be captured, so that such a problem need only be solved under whichever of the two ambiguity sets is more convenient.

\section*{Acknowledgments}
The claims and proofs for mean--MAD ambiguity, including Theorem~\ref{th:bijection}, were obtained by the authors without the use of generative AI. For the mean--variance problem, we initially had an exact characterization only for $N=2$, which already required a long technical proof, and we did not see how to extend this to arbitrary $N\in \N$. Around June 2026, Claude Opus produced a complete proof for arbitrary $N$ (now Theorem~\ref{th:var_N}). However, to us its exposition was so exceedingly poor that it took us considerable time to understand its main arguments, and considerably more to simplify and rewrite. The current proof is our own exposition, and is substantially different from the original AI-generated proof. The proofs of some of the asymptotic results under mean--variance ambiguity were generated with AI and checked and edited by the authors. Finally, Claude Opus 5.5 and ChatGPT-6 Astra were used for language editing and for coding the numerical illustrations, and all code and numerical results were checked by the authors. The authors take full responsibility for all content of this paper.

\begin{small}
\bibliography{bibbook}
\end{small}

\begin{appendix}
\section{Proof of Lemma \ref{lem:tail_equivalence}}\label{ap:Lemma_taileq}
Since $\bP(X>p_i)\le\bP(X\ge p_i)$ for every $\bP\in\cP_\sigma(\mu,\sigma)$ and $i\in\{1,\dots,N\}$, it holds that
\[
    \sup_{(w,p)\in\mathcal H}\ \inf_{\bP\in\cP_\sigma(\mu,\sigma)}\sum_{i=1}^Nw_ip_i\bP(X>p_i)\ \le\ \sup_{(w,p)\in\mathcal H}\ \inf_{\bP\in\cP_\sigma(\mu,\sigma)}\sum_{i=1}^Nw_ip_i\bP(X\ge p_i).
\]
Conversely, fix $(w,p)\in\mathcal H$ and $\varepsilon\in(0,1)$. Since $(w,(1-\varepsilon)p)\in\mathcal H$ and $\bP(X>(1-\varepsilon)p_i)\ge\bP(X\ge p_i)$,
\begin{align*}
    \sup_{(w,p)\in\mathcal H}\ \inf_{\bP\in\cP_\sigma(\mu,\sigma)}\sum_{i=1}^Nw_ip_i\bP(X>p_i)
    &\ \ge\ \inf_{\bP\in\cP_\sigma(\mu,\sigma)}\sum_{i=1}^Nw_i(1-\varepsilon)p_i\,\bP\bigl(X>(1-\varepsilon)p_i\bigr)\\
    &\ \ge\ (1-\varepsilon)\inf_{\bP\in\cP_\sigma(\mu,\sigma)}\sum_{i=1}^Nw_ip_i\bP(X\ge p_i).
\end{align*}
Taking the supremum over $(w,p)\in\mathcal H$ and letting $\varepsilon\downarrow0$ gives
\[
    \sup_{(w,p)\in\mathcal H}\ \inf_{\bP\in\cP_\sigma(\mu,\sigma)}\sum_{i=1}^Nw_ip_i\bP(X>p_i)\ \ge\ \sup_{(w,p)\in\mathcal H}\ \inf_{\bP\in\cP_\sigma(\mu,\sigma)}\sum_{i=1}^Nw_ip_i\bP(X\ge p_i).
\]

Moreover, if $(w,p)\in\mathcal H$ is optimal under the strict-tail convention, then
\[
    \inf_{\bP\in\cP_\sigma(\mu,\sigma)}\sum_{i=1}^Nw_ip_i\bP(X\ge p_i)\ \ge\ \inf_{\bP\in\cP_\sigma(\mu,\sigma)}\sum_{i=1}^Nw_ip_i\bP(X>p_i),
\]
so that $(w,p)$ attains the common value under the weak-tail convention as well.

\section{Proof of Lemma \ref{lem:e_and_u}}\label{ap:Lemma_var}
Fix $N\ge2$ and let $\mathcal S_p$ denote the set of pairs $(p_1,p_2)\in\mathbb{R}^2$ for which \eqref{eq:zdef} and \eqref{eq:rec} generate quantities $(r_1,r_2,p_3,r_3,\dots,p_N,r_N,z)$ such that $0<p^*_1<\dots<p^*_N$, and such that \eqref{eq:terminal} and \eqref{eq:moment} hold. The outline of the proof is as follows. First, we show that the map $S_p\mapsto p_N/z$ is injective with $p_N/z \in(0,1)$. Second, we show that every such value solves an equation $\Phi(x)=\mu^2/(\mu^2+\sigma^2)$. Finally, we show that $\Phi$ is strictly increasing on $(0,1)$, so that the equation has at most one solution in $(0,1)$. Hence, injectivity implies that $\mathcal S_p$ has at most one element.
 
Fix a pair $(p_1,p_2)\in\mathcal S_p$ and let $r_1,r_2,p_3,r_3,\dots,p_N,r_n,z$ be its generated quantities. Squaring \eqref{eq:rec} gives
\begin{align*}
    (z-p_{j+1})^2=(z-p_j)(z+p_j-2r_{j-1}p_j)
\end{align*}
for $j=2,\dots,N-1$, and by \eqref{eq:terminal} the same identity holds at $j=N$, where both sides are zero because $p_{N+1}=z$. Hence, it holds for all $j = 2,\dots,N$ that
\begin{align}\label{eq:claim1}
    (z-p_{j+1})^2=(z-p_j)(z+p_j-2r_{j-1}p_j).
\end{align}

We now normalize by $z$, so that $p_{N+1}/z=1$ and $r_i=(p_{i+1}/z)/(p_i/z)$. Dividing \eqref{eq:claim1} by $z^2$ gives $(1-p_{j+1}/z)^2=1-(p_j/z)^2-2(p_j/z)^2(1-p_j/z)/(p_{j-1}/z)$, and since $1-(1-p_{j+1}/z)^2=(p_{j+1}/z)(2-p_{j+1}/z)$, it holds that
\begin{align}\label{eq:backward}
    \frac{p_{j-1}}{z}=\frac{2(p_j/z)^2(1-p_j/z)}{(p_{j+1}/z)(2-p_{j+1}/z)-(p_j/z)^2}.
\end{align} Moreover, \eqref{eq:terminal} gives $p_N/z=1/(2r_{N-1}-1)$, so that $p_N/z\in(0,1)$.

To compare different solutions, we also evaluate \eqref{eq:backward} from an arbitrary starting value. For $x\in(0,1)$, set $\rho_{N+1}(x):=1$ and $\rho_N(x):=x$, define $\rho_{N-1}(x),\dots,\rho_1(x)$ by \eqref{eq:backward} with $\rho_j(x)$ in place of $p_j/z$, and let $\rho_0(x)$ be the auxiliary quantity obtained in this way at $j=1$. The normalized prices of the solution are $p_j/z=\rho_j(p_N/z)$ for $j=1,\dots,N+1$. Write $s_j(x):=\rho_j(x)\bigl(2-\rho_j(x)\bigr)$. If $0<\rho_j(x)<\rho_{j+1}(x)\le1$, then $s_{j+1}(x)>s_j(x)$, since the function $t\mapsto t(2-t)$ is increasing on $(0,1]$, so that $s_{j+1}(x)-\rho_j(x)^2>s_j(x)-\rho_j(x)^2=2\rho_j(x)\bigl(1-\rho_j(x)\bigr)>0$ and $0<\rho_{j-1}(x)<\rho_j(x)$. Backward induction from $\rho_N(x)=x<1=\rho_{N+1}(x)$ therefore gives $0<\rho_0(x)<\rho_1(x)<\dots<\rho_N(x)<\rho_{N+1}(x)=1$, and since the denominators are all positive, each $\rho_j(x)$ is differentiable with respect to $x$.

For $x\in(0,1)$, let $r_i(x):=\rho_{i+1}(x)/\rho_i(x)$,
\[
    S(x):=\sum_{i=1}^N\bigl(r_i(x)-1\bigr),\qquad T(x):=\sum_{i=1}^N\rho_i(x)\bigl(r_i(x)^2-1\bigr),
\]
and notice that $T(x)>0$. Furthermore, define
\[
    M(x):=\rho_1(x)+\frac{2\rho_1(x)\bigl(1-\rho_1(x)\bigr)S(x)}{T(x)},\qquad \Phi(x):=\frac{M(x)^2}{s_1(x)}.
\]

We now return to the fixed pair $(p_1,p_2)\in\mathcal S_p$. So far we have used \eqref{eq:rec} and \eqref{eq:terminal}. The remaining conditions, \eqref{eq:zdef} and \eqref{eq:moment}, give two equations in $p_N/z$ and $z$. Since \eqref{eq:zdef} states that $\mu^2+\sigma^2=2zp_1-p_1^2$, dividing by $z^2$ and using $p_1/z=\rho_1(p_N/z)$ gives $\frac{\mu^2+\sigma^2}{z^2}=s_1(p_N/z)$. Similarly, since $p_j/z=\rho_j(p_N/z)$ for all $j$, dividing \eqref{eq:moment} by $z^2$ gives $$\bigl(s_1(p_N/z)-\rho_1(p_N/z)^2\bigr)S(p_N/z)=\Bigl(\frac{\mu}{z}-\rho_1(p_N/z)\Bigr)T(p_N/z).$$ In conclusion, it holds that 
\begin{align}\label{eq:norm1}
    \frac{\mu^2+\sigma^2}{z^2}=s_1(p_N/z)
\end{align}
and that
\begin{align}\label{eq:norm2}
    \bigl(s_1(p_N/z)-\rho_1(p_N/z)^2\bigr)S(p_N/z)=\Bigl(\frac{\mu}{z}-\rho_1(p_N/z)\Bigr)T(p_N/z).
\end{align}

Since $s_1(p_N/z)-\rho_1(p_N/z)^2=2\rho_1(p_N/z)\bigl(1-\rho_1(p_N/z)\bigr)$ and $T(p_N/z)>0$, solving \eqref{eq:norm2} for $\mu/z$ gives
\begin{align}\label{eq:muz}
    \frac{\mu}{z}=\rho_1(p_N/z)+\frac{2\rho_1(p_N/z)\bigl(1-\rho_1(p_N/z)\bigr)S(p_N/z)}{T(p_N/z)}=M(p_N/z).
\end{align}
Together with \eqref{eq:norm1}, this gives
\begin{align}\label{eq:norm3}
    \Phi(p_N/z)=\frac{M(p_N/z)^2}{s_1(p_N/z)}=\frac{(\mu/z)^2}{(\mu^2+\sigma^2)/z^2}=\frac{\mu^2}{\mu^2+\sigma^2}.
\end{align}

It remains to show that $\Phi$ is strictly increasing on $(0,1)$. To establish injectivity, take two pairs $(p_1,p_2),(\tilde p_1,\tilde p_2)\in\mathcal S_p$, with associated generated quantities $z,p_3,\dots,p_N$ and $\tilde z,\tilde p_3,\dots,\tilde p_N$, and suppose that
\[
    \frac{p_N}{z}=\frac{\tilde p_N}{\tilde z}=x.
\]
Then $p_j/z=\rho_j(x)=\tilde p_j/\tilde z$ for $j=1,\dots,N+1$, and the identity $\mu/z=M(x)$ above, applied to both pairs, gives
\[
    z=\frac{\mu}{M(x)}=\tilde z.
\]
Therefore $p_j=\tilde p_j$ for $j=1,\dots,N$, and in particular $(p_1,p_2)=(\tilde p_1,\tilde p_2)$. Thus the map is injective on $\mathcal S_p$. By \eqref{eq:norm3}, the image of every pair in $\mathcal S_p$ solves $\Phi(x)=\mu^2/(\mu^2+\sigma^2)$. If $\Phi$ is strictly increasing on $(0,1)$, this equation has at most one solution there, so all pairs in $\mathcal S_p$ have the same image and injectivity implies that $\mathcal S_p$ has at most one element.

Since $\Phi$ is a differentiable function on $(0,1)$, it remains to show that $\Phi'(x)>0$ for all $x\in(0,1)$. Since $\rho_1(x)\in(0,1)$ and $S(x),T(x)>0$, both $M(x)>0$ and $s_1(x)>0$. Hence, differentiating $\Phi(x)$ results in
\begin{align}\label{eq:claim8}
    \frac{\Phi'(x)}{\Phi(x)}=\frac{\bigl(2M'(x)-s_1'(x)\bigr)-\dfrac{s_1'(x)}{s_1(x)}\bigl(M(x)-s_1(x)\bigr)}{M(x)}.
\end{align}
Since $\Phi(x)>0$ and $M(x)>0$, it follows that $\Phi'(x)>0$ if and only if the numerator in \eqref{eq:claim8} is strictly positive. This numerator involves $M(x)-s_1(x)$ and $2M'(x)-s_1'(x)$, and we start with the former.

Recall that $M(x)=\rho_1(x)+2\rho_1(x)\bigl(1-\rho_1(x)\bigr)S(x)/T(x)$. To simplify $M(x)-s_1(x)$, we first compute $T(x)$ through its partial sums. For $j=N$, since $\rho_{N+1}(x)=s_{N+1}(x)=1$, we have $\rho_N(x)\bigl(r_N(x)^2-1\bigr)=\bigl(s_{N+1}(x)-\rho_N(x)^2\bigr)/\rho_N(x)$. Moreover, if $\sum_{i=j}^N\rho_i(x)\bigl(r_i(x)^2-1\bigr)=\bigl(s_{j+1}(x)-\rho_j(x)^2\bigr)/\rho_j(x)$ for some $j\ge2$, then by \eqref{eq:backward} for $\rho_j(x)$ this sum equals $2\rho_j(x)\bigl(1-\rho_j(x)\bigr)/\rho_{j-1}(x)$, and adding $\rho_{j-1}(x)\bigl(r_{j-1}(x)^2-1\bigr)=\bigl(\rho_j(x)^2-\rho_{j-1}(x)^2\bigr)/\rho_{j-1}(x)$ gives the same identity at $j-1$. Hence, by backward induction,
\begin{align}\label{eq:norm4}
    \sum_{i=j}^N\rho_i(x)\bigl(r_i(x)^2-1\bigr)=\frac{s_{j+1}(x)-\rho_j(x)^2}{\rho_j(x)},
\end{align}
for all $j = 1,\dots,N$. At $j=1$, the sum in \eqref{eq:norm4} is $T(x)$, and \eqref{eq:backward} for $\rho_j(x)$ at $j=1$ gives $s_2(x)-\rho_1(x)^2=2\rho_1(x)^2\bigl(1-\rho_1(x)\bigr)/\rho_0(x)$, so that
\[
    T(x)=\frac{s_2(x)-\rho_1(x)^2}{\rho_1(x)}=\frac{2\rho_1(x)\bigl(1-\rho_1(x)\bigr)}{\rho_0(x)}.
\]
Therefore,
\begin{align}\label{eq:norm5}
    \rho_0(x)T(x)=2\rho_1(x)\bigl(1-\rho_1(x)\bigr).
\end{align}

Substituting \eqref{eq:norm5} into the definition of $M(x)$ gives \begin{align}\label{eq:Mformula}
    M(x)=\rho_1(x)+\rho_0(x)S(x),
\end{align} and since $s_1(x)=\rho_1(x)\bigl(2-\rho_1(x)\bigr)$, \eqref{eq:norm5} gives
\[
    M(x)-s_1(x)=\rho_0(x)S(x)-\rho_1(x)\bigl(1-\rho_1(x)\bigr)=\rho_0(x)S(x)-\tfrac12\rho_0(x)T(x) = \tfrac12\rho_0(x)\bigl(2S(x)-T(x)\bigr).
\]
Hence, we obtain
\begin{align}\label{eq:MD}
    M(x)-s_1(x)=\tfrac12\rho_0(x)\bigl(2S(x)-T(x)\bigr).
\end{align}

To compute $2M'(x)-s_1'(x)$, we differentiate \eqref{eq:MD}. Writing $D(x):=2S(x)-T(x)$, this gives
\[
    M'(x)-s_1'(x)=\tfrac12\rho_0'(x)D(x)+\tfrac12\rho_0(x)D'(x),
\]
so we first compute $D'(x)$. Since $\bigl(r_i(x)-1\bigr)\bigl(\rho_i(x)+\rho_{i+1}(x)\bigr)=\rho_i(x)\bigl(r_i(x)^2-1\bigr)$, we can write
\begin{align}\label{eq:norm6}
    D(x)=\sum_{i=1}^N\frac{\bigl(1-\rho_i(x)\bigr)^2-\bigl(1-\rho_{i+1}(x)\bigr)^2}{\rho_i(x)}.
\end{align}

By \eqref{eq:norm6}, $\rho_j(x)$ with $2\le j\le N$ appears only in the terms $i=j-1$ and $i=j$, $\rho_1(x)$ only in the term $i=1$, and $\rho_{N+1}(x)=1$ is constant. Differentiating and collecting terms therefore gives
\begin{align}\label{eq:norm7}
    D'(x)=\sum_{j=1}^Nc_j(x)\rho_j'(x),
\end{align}
 where
\[
    c_1(x)=-\frac{s_2(x)-\rho_1(x)^2}{\rho_1(x)^2},\qquad
    c_j(x)=\frac{2\bigl(1-\rho_j(x)\bigr)}{\rho_{j-1}(x)}-\frac{s_{j+1}(x)-\rho_j(x)^2}{\rho_j(x)^2},\quad 2\le j\le N.
\]

By \eqref{eq:backward} for $\rho_j(x)$, $c_j(x)=0$ for $2\le j\le N$, and $c_1(x)=-2\bigl(1-\rho_1(x)\bigr)/\rho_0(x)$ by \eqref{eq:backward} at $j=1$. Hence, $D'(x)=-2\bigl(1-\rho_1(x)\bigr)\rho_1'(x)/\rho_0(x)$, and since $s_1'(x)=2\bigl(1-\rho_1(x)\bigr)\rho_1'(x)$, it holds that $\rho_0(x)D'(x)=-s_1'(x)$.

Differentiating $M(x)-s_1(x)=\tfrac12\rho_0(x)D(x)$ gives
\[
    M'(x)-s_1'(x)=\tfrac12\rho_0'(x)D(x)+\tfrac12\rho_0(x)D'(x)=\tfrac12\rho_0'(x)D(x)-\tfrac12s_1'(x),
\]
so that $2M'(x)-s_1'(x)=\rho_0'(x)D(x)$. With $\ell_j(x):=\rho_j'(x)/\rho_j(x)$, we have $\rho_0'(x)=\ell_0(x)\rho_0(x)$, and this is equivalent to $\rho_0(x)D(x)=2\bigl(M(x)-s_1(x)\bigr)$. Hence, we obtain
\begin{align}\label{eq:claim15}
    2M'(x)-s_1'(x)=2\ell_0(x)\bigl(M(x)-s_1(x)\bigr).
\end{align}

Substituting \eqref{eq:claim15} into the numerator in \eqref{eq:claim8}, and using $M(x)-s_1(x)=\tfrac12\rho_0(x)D(x)$ and $s_1'(x)/s_1(x)=2\bigl(1-\rho_1(x)\bigr)\ell_1(x)/\bigl(2-\rho_1(x)\bigr)$, gives 
\begin{align}\label{eq:claim16}
    \frac{\Phi'(x)}{\Phi(x)}=\frac{\rho_0(x)D(x)}{M(x)}\Bigl(\ell_0(x)-\frac{1-\rho_1(x)}{2-\rho_1(x)}\,\ell_1(x)\Bigr).
\end{align}

In \eqref{eq:claim16}, the first factor is positive, since $\rho_0(x)>0$, $D(x) > 0$, and $M(x)>0$. It remains to show that $\ell_0(x)>\ell_1(x)>0$, because then
\[
    \ell_0(x)-\frac{1-\rho_1(x)}{2-\rho_1(x)}\,\ell_1(x)>\ell_0(x)-\ell_1(x)>0,
\]
as $\bigl(1-\rho_1(x)\bigr)/\bigl(2-\rho_1(x)\bigr)<1$.

We have $\ell_{N+1}(x)=0$ and $\ell_N(x)=1/x$. All factors in \eqref{eq:backward} for $\rho_j(x)$ are positive, and taking logarithms gives
\[
    \ln\rho_{j-1}(x)=\ln2+2\ln\rho_j(x)+\ln\bigl(1-\rho_j(x)\bigr)-\ln\bigl(s_{j+1}(x)-\rho_j(x)^2\bigr),\qquad j=1,\dots,N.
\]
Differentiating and rearranging gives for each $j=1,\dots,N$ that
\[
    \ell_{j-1}(x)-\ell_j(x)=a_j(x)\bigl(\ell_j(x)-\ell_{j+1}(x)\bigr)+b_j(x)\ell_j(x),
\]
where
\[
    a_j(x):=\frac{2\rho_{j+1}(x)\bigl(1-\rho_{j+1}(x)\bigr)}{s_{j+1}(x)-\rho_j(x)^2},\qquad
    b_j(x):=\frac{\bigl(\rho_{j+1}(x)-\rho_j(x)\bigr)^2}{\bigl(s_{j+1}(x)-\rho_j(x)^2\bigr)\bigl(1-\rho_j(x)\bigr)}.
\]

Clearly, $a_j(x)\ge0$ and $b_j(x)>0$. Therefore, $\ell_j(x)>\ell_{j+1}(x)\ge0$ implies $\ell_{j-1}(x)>\ell_j(x)$, and backward induction from $\ell_N(x)>\ell_{N+1}(x)=0$ shows that $0<\ell_N(x)<\dots<\ell_1(x)<\ell_0(x)$. This proves that $\Phi'(x)>0$ for all $x\in(0,1)$, so $\Phi$ is strictly increasing on $(0,1)$ and $\mathcal S_p$ has at most one element. This completes the proof.

\section{Proof of Proposition \ref{prop:prices_var}}\label{ap:prices_var}
Recall from the proof of Theorem~\ref{th:var_N} and Lemma~\ref{lem:tail_equivalence} that $R^*=\Psi^*=\max_{p,z}\Psi(p,z)$, with $\Psi(p,z)=A(p,z)/D(p,z)$ as in \eqref{eq:Psi}. When we wish to make the dependence on $\sigma$ explicit, we write it as an additional argument, e.g., $\Psi(p,z;\sigma)$ or $R^*(\sigma)$.

We start with Part~(ii). The domain of $\Psi$ and the denominator $D(p,z)>0$ do not depend on $\sigma$, whereas $A(p,z;\sigma)=(z-p_1)^2-(z-\mu)^2-\sigma^2$ is strictly decreasing in $\sigma$, so that $\Psi(p,z;\sigma)$ is strictly decreasing in $\sigma$ as well. Let $\sigma_1<\sigma_2$, and let $(p,z)$ attain the maximum of $\Psi(\cdot\,;\sigma_2)$. Then
\[
    R^*(\sigma_2)=\Psi(p,z;\sigma_2)<\Psi(p,z;\sigma_1)\le R^*(\sigma_1),
\]
so that $R^*$ is strictly decreasing in $\sigma$.

We now prove Part~(i). By \eqref{eq:weights_var}, $w^*_ip^*_i=\bigl((z^*-p^*_i)^2-(z^*-p^*_{i+1})^2\bigr)/D(p^*,z^*)$, so that summing over $i$ gives
\begin{align}\label{eq:EP_formula}
    \bE(P;\sigma)=\sum_{i=1}^Nw^*_ip^*_i=\frac{(z^*-p^*_1)^2}{D(p^*,z^*)}.
\end{align}

We first show that $\lim_{\sigma\downarrow0}\bE(P;\sigma)=\mu$ with $\bE(P;\sigma)<\mu$ for all sufficiently small $\sigma$. Since $\bE(P;\sigma)$ is a weighted average of the prices, $p^*_1\le\bE(P;\sigma)\le p^*_N$, so that it suffices to show that $p^*_1\to\mu$ as $\sigma\downarrow0$, and that $p^*_N<\mu$ for all sufficiently small $\sigma$.

We start with $p^*_1$. Since $w^*_1>0$, \eqref{eq:intervalmasses} implies that $m_1=1-R^*/p^*_1\ge0$, so that $R^*\le p^*_1$. Moreover, $p^*_1<\mu$ by \eqref{eq:rev_var_N}, so that $R^*\le p^*_1<\mu$. Next, for $\sigma<\mu$, Cantelli's inequality shows that the posted price $p=\mu-\sqrt{\mu\sigma}$ guarantees
\[
    \frac{p(\mu-p)^2}{(\mu-p)^2+\sigma^2}=\frac{p\,\mu\sigma}{\mu\sigma+\sigma^2}=\frac{(\mu-\sqrt{\mu\sigma})\,\mu}{\mu+\sigma},
\]
which converges to $\mu$, so that indeed $R^*\to\mu$. Hence, $p^*_1\to\mu$, and then also $z^*\to\mu$ by \eqref{eq:zdef}.

We now show that $p^*_N<\mu$ for all sufficiently small $\sigma$. Recall from \eqref{eq:sqrd_eq} that
\[
    (z^*-p^*_{i+1})^2=(z^*-p^*_i)\bigl(z^*+p^*_i-2r^*_{i-1}p^*_i\bigr),\qquad i=2,\dots,N.
\]
The left-hand side is non-negative, and $z^*-p^*_i>0$, so that the second factor on the right-hand side is non-negative as well, i.e.,
\[
    z^*-p^*_i\ \ge\ 2r^*_{i-1}p^*_i-2p^*_i=2\,\frac{p^*_i}{p^*_{i-1}}\,(p^*_i-p^*_{i-1})\ >\ 2(p^*_i-p^*_{i-1}),
\]
where the last inequality uses $p^*_i>p^*_{i-1}$. Hence,
\[
    z^*-p^*_{i-1}=(z^*-p^*_i)+(p^*_i-p^*_{i-1})\ <\ (z^*-p^*_i)+\tfrac12(z^*-p^*_i)=\tfrac32(z^*-p^*_i),\qquad i=2,\dots,N.
\]
Applying this repeatedly, starting from $i=2$, gives
\[
    z^*-p^*_1\ <\ \tfrac32(z^*-p^*_2)\ <\ \bigl(\tfrac32\bigr)^2(z^*-p^*_3)\ <\ \dots\ <\ \bigl(\tfrac32\bigr)^{N-1}(z^*-p^*_N),
\]
so that
\[
    z^*-p^*_N\ >\ \bigl(\tfrac23\bigr)^{N-1}(z^*-p^*_1).
\]
On the other hand, multiplying \eqref{eq:zdef} by $2p^*_1$ gives $2z^*p^*_1-(p^*_1)^2=\mu^2+\sigma^2$, so that
\[
    A(p^*,z^*)=2\mu z^*-\mu^2-\sigma^2-\bigl(2z^*p^*_1-(p^*_1)^2\bigr)=2(\mu z^*-\mu^2-\sigma^2).
\]
Since $A(p^*,z^*)=R^*D(p^*,z^*)>0$, it follows that $\mu z^*>\mu^2+\sigma^2=(z^*)^2-(z^*-p^*_1)^2$, and hence $0<z^*-\mu<(z^*-p^*_1)^2/z^*$. As $(z^*-p^*_1)/z^*\to0$, we have $(z^*-p^*_1)/z^*<(2/3)^{N-1}$ for all sufficiently small $\sigma$, and then
\[
    z^*-\mu\ <\ \bigl(\tfrac23\bigr)^{N-1}(z^*-p^*_1)\ <\ z^*-p^*_N,
\]
which implies that $p^*_N<\mu$.

Next, we show that $\lim_{\sigma\to\infty}\bE(P;\sigma)=\infty$, by showing that $\bE(P;\sigma)\ge p^*_N/(2N)$ and $p^*_N>\tfrac12\sqrt{\mu^2+\sigma^2}$. Substituting $p^*_{i+1}=r^*_ip^*_i$ into \eqref{eq:sqrd_eq} and rearranging gives
\[
    2(z^*-p^*_i)(r^*_i-r^*_{i-1})=p^*_i(r^*_i-1)^2>0,\qquad i=2,\dots,N,
\]
so that $r^*_1<\dots<r^*_N$. Since $2z^*-p^*_i-p^*_{i+1}\le2(z^*-p^*_1)$ and $r^*_N-1=(z^*-p^*_N)/p^*_N\le(z^*-p^*_1)/p^*_N$, every term of $D(p^*,z^*)$ is at most $2(z^*-p^*_1)^2/p^*_N$, so that \eqref{eq:EP_formula} gives $\bE(P;\sigma)\ge p^*_N/(2N)$. Moreover, \eqref{eq:terminal} gives $z^*+p^*_N=2(p^*_N)^2/p^*_{N-1}$, so that, as $p^*_{N-1}\ge p^*_1$ and $2z^*p^*_1=(p^*_1)^2+\mu^2+\sigma^2$,
\[
    (p^*_N)^2=\tfrac12p^*_{N-1}(z^*+p^*_N)\ >\ \tfrac12p^*_1z^*\ >\ \tfrac14(\mu^2+\sigma^2).
\]
Hence, $\bE(P;\sigma)>\sqrt{\mu^2+\sigma^2}/(4N)$, which diverges as $\sigma\to\infty$.

It remains to show that the derivative of $\bE(P;\sigma)$ with respect to $\sigma$ has only finitely many zeros, so that $\bE(P;\sigma)$ is strictly monotone on $(0,\sigma_l)$ and on $(\sigma_h,\infty)$ for some $0<\sigma_l\le\sigma_h<\infty$. Since $\bE(P;\sigma)$ approaches $\mu$ from below as $\sigma\downarrow0$ and diverges as $\sigma\to\infty$, it is then decreasing on the first interval and increasing on the second.

The key observation is that $\bE(P;\sigma)$ depends on $\sigma$ only through the ratio $g(\sigma):=p^*_N/z^*$, i.e., $\bE(P;\sigma)=\mu f(g(\sigma))$ for a rational function $f$ that does not depend on $\sigma$. We show that $f'$ has only finitely many zeros and that $g'(\sigma)<0$ for all $\sigma>0$, so that the claim follows from the chain rule.

We start with $f$. Recall from the proof of Lemma~\ref{lem:e_and_u} that $p^*_j/z^*=\rho_j(g(\sigma))$ for $j=1,\dots,N+1$, where $\rho_{N+1}(x)=1$, $\rho_N(x)=x$ and
\[
    \rho_{j-1}(x)=\frac{2\rho_j(x)^2\bigl(1-\rho_j(x)\bigr)}{\rho_{j+1}(x)\bigl(2-\rho_{j+1}(x)\bigr)-\rho_j(x)^2},\qquad j=N,\dots,2,
\]
and that $\mu/z^*=M(g(\sigma))$ by \eqref{eq:muz}. Together with \eqref{eq:EP_formula}, this gives $\bE(P;\sigma)=\mu f(g(\sigma))$ with
\[
    f(x):=\frac{\bigl(1-\rho_1(x)\bigr)^2}{M(x)\sum_{i=1}^N\bigl(\frac{\rho_{i+1}(x)}{\rho_i(x)}-1\bigr)\bigl(2-\rho_i(x)-\rho_{i+1}(x)\bigr)}.
\]
Since $\rho_1,\dots,\rho_{N+1}$ and $M$ are rational functions that do not depend on $\sigma$, so is $f$, and its denominator is positive on $(0,1)$. Moreover, $f$ is not constant, because $f(g(\sigma))=\bE(P;\sigma)/\mu$ converges to $1$ as $\sigma\downarrow0$ and diverges as $\sigma\to\infty$. Hence, $f'$ has only finitely many zeros in $(0,1)$.

We now turn to $g$. Recall from \eqref{eq:norm3} that
\[
    \Phi\bigl(g(\sigma)\bigr)=\frac{\mu^2}{\mu^2+\sigma^2},
\]
where $\Phi$ is differentiable with $\Phi'(x)>0$ for all $x\in(0,1)$. Hence, $\Phi$ has a differentiable inverse, so that $g$ is differentiable as well, and differentiating both sides with respect to $\sigma$ gives
\[
    \Phi'\bigl(g(\sigma)\bigr)\,g'(\sigma)=-\frac{2\mu^2\sigma}{(\mu^2+\sigma^2)^2}.
\]
Since the right-hand side is negative and $\Phi'(g(\sigma))>0$, it follows that $g'(\sigma)<0$.

Finally, the chain rule gives
\[
    \frac{d\,\bE(P;\sigma)}{d\sigma}=\mu\,f'\bigl(g(\sigma)\bigr)\,g'(\sigma),
\]
which is zero if and only if $f'(g(\sigma))=0$, as $g'(\sigma)<0$. Since $g$ is strictly decreasing, it attains each of the finitely many zeros of $f'$ for at most one $\sigma$, so that the derivative of $\bE(P;\sigma)$ has only finitely many zeros. This finishes the proof.

\section{Proof of Proposition \ref{prop:prices_decrease_MAD}}\label{ap:monotonicity}
Recall that $p_1(N)$ is the unique solution on the interval $(0,\mu-d/2)$ to
    \begin{align}\label{eq:p_1eq}
        \frac{\mu-d/2}{p_1(N)}+N\left(\frac{p_1(N)}{\mu}\right)^{1/N}=N+1.
    \end{align}
Then, by writing $k = \mu-d/2$ and differentiating \eqref{eq:p_1eq} with respect to $k$, we obtain
\begin{align}
    \frac{1}{p_1(N)} + \left(-\frac{k}{p_1(N)^2}+\frac{1}{p_1(N)}\left(\frac{p_1(N)}{\mu}\right)^{\frac{1}{N}}\right)\frac{\partial}{\partial k}p_1(N) = 0.
\end{align}
Hence,
\begin{align}
    \frac{\partial}{\partial k}p_1(N) = \frac{p_1(N)}{k - p_1(N)(p_1(N)/\mu)^{1/N}}.
\end{align}
Then, since $p_1(N)(p_1(N)/\mu)^{1/N} \leq p_1(N) < k$ and $p_1(N)>0$, we have that $\frac{\partial}{\partial k}p_1(N) > 0$. Therefore, since $k$ is decreasing in $d$, we have that $p_1(N)$ is decreasing in $d$.

For the other prices, we have
    \begin{align}
        p_j(N) = p_1(N)\left(\frac{\mu}{p_1(N)}\right)^{\frac{j-1}{N}} = \mu^{\frac{j-1}{N}}p_1(N)^{\frac{N-j+1}{N}},
        \quad j=2,\dots,N.
    \end{align}
Since $(N-j+1)/N > 0$, we have that $p_j(N)$ is increasing in $p_1(N)$, which in turn makes $p_j(N)$ decreasing in $d$. 

Since the weights are uniform, the expected price $\bE(P;d)=\frac1N\sum_{i=1}^Np_i(N)$ is decreasing in $d$ as well, which proves Part~(i). Part~(ii) follows from Part~(i), as $R^*=p_1(N)p_N(N)/\mu$ by \eqref{eq:rev_N_2a}.

\section{Asymptotic properties of the finite menu}\label{ap:asympt}
This appendix shows that the optimal finite menu recovers full randomization as the menu size grows, first under mean--variance ambiguity and then under mean--MAD ambiguity. In both cases, we prove that the lowest optimal price $p^*_1(N)$ and the optimal expected revenue $R^*(N)$ converge to the optimal expected revenue of the fully randomized mechanism, and that the selling mechanism
\begin{align}\label{eq:seller_dist}
    \bP^*_{\mathrm{seller}}(N)=\sum_{i=1}^Nw^*_i\,\delta_{p^*_i(N)}
\end{align}
converges in distribution, where $(w^*,p^*)$ is the optimal menu of Theorem~\ref{th:var_N} or Theorem~\ref{th:mad_N}, respectively, and $\delta_x$ is the Dirac measure at $x$. Throughout, we sometimes write $R^*(N)$, $p^*_i(N)$, and $z^*(N)$ to highlight the dependence on the menu size. Moreover, $\xrightarrow{d}$ denotes convergence in distribution.

\subsection{Mean--variance convergence}
The next proposition establishes this convergence for the mean--variance ambiguity set $\cP_\sigma(\mu,\sigma)$, for which the fully randomized revenue $\pi_\sigma$ in \eqref{eq:pisigma} follows from \citet{carrasco2018optimal}, together with the convergence of the selling mechanism \eqref{eq:seller_dist}.

\begin{proposition}\label{prop:asympt_var}
    Consider $\cP=\cP_\sigma(\mu,\sigma)$ and define $z_\infty:=\pi_\sigma e^{\mu/\pi_\sigma-1}$. The following asymptotic properties hold:
    \begin{enumerate}[label=\textup{(\roman*)}]
        \item $\lim_{N\to\infty}p^*_1(N)=\pi_\sigma$;
        \item $\lim_{N\to\infty}R^*(N)=\pi_\sigma$;
        \item $\bP^*_{\mathrm{seller}}(N) \xrightarrow{d} \mathbb{H}^*_{\mathrm{seller}}$, where $\mathbb{H}^*_{\mathrm{seller}}$ has CDF
        \begin{align}
            H^*(p) =
            \begin{cases}
                0, & p<\pi_\sigma, \\[1mm]
                \dfrac{z_\infty\ln(p/\pi_\sigma)-(p-\pi_\sigma)}{z_\infty\ln(z_\infty/\pi_\sigma)-(z_\infty-\pi_\sigma)}, & p\in[\pi_\sigma,z_\infty], \\[2mm]
                1, & p>z_\infty.
            \end{cases}
        \end{align}
    \end{enumerate}
\end{proposition}

\begin{proof}
Fix $N\ge2$. We first show that
\begin{align}\label{eq:p1_gap}
    0<p^*_1-R^*<\frac{\mu-p^*_1}{N},
\end{align}
so that $p^*_1(N)$ and $R^*(N)$ have the same limit, provided that $R^*(N)$ converges. We then show that $R^*(N)$ indeed converges, and that its limit solves \eqref{eq:pisigma}. Finally, we show that the selling mechanism converges.

Recall from the proof of Theorem~\ref{th:var_N} that $(p^*,z^*)$ maximizes $\Psi(p,z)=A(p,z)/D(p,z)$ with $\Psi(p^*,z^*)=R^*$, where
\[
    \frac{\partial A(p,z)}{\partial p_1}=-2(z-p_1)\qquad\text{and}\qquad\frac{\partial D(p,z)}{\partial p_1}=1-\frac{2zp_2-p_2^2}{p_1^2}.
\]
The first-order condition with respect to $p_1$ therefore reads
\[
    2(z^*-p^*_1)=R^*\biggl(\frac{2z^*p^*_2-(p^*_2)^2}{(p^*_1)^2}-1\biggr).
\]
Since $2z^*p^*_2-(p^*_2)^2-(p^*_1)^2=(p^*_2-p^*_1)(2z^*-p^*_1-p^*_2)+2p^*_1(z^*-p^*_1)$, rearranging gives
\[
    \frac{p^*_1}{R^*}-1=(r^*_1-1)\,\frac{2z^*-p^*_1-p^*_2}{2(z^*-p^*_1)}.
\]
The right-hand side is positive, and since $p^*_2>p^*_1$, it is smaller than $r^*_1-1$. Moreover, recall from the proof of Proposition~\ref{prop:prices_var} that $r^*_1<\dots<r^*_N$, so that $r^*_1-1$ is at most the average of $r^*_1-1,\dots,r^*_N-1$, which equals $(\mu-p^*_1)/(NR^*)$ by \eqref{eq:rev_var_N}. This proves \eqref{eq:p1_gap}.

Next, $R^*(N)$ is non-decreasing in $N$, since a menu with $N$ prices can be written as a menu with $N+1$ prices by repeating a price, and $R^*(N)<p^*_1(N)<\mu$ by \eqref{eq:p1_gap} and \eqref{eq:rev_var_N}. Hence, $R^*(N)$ converges to some $\ell\in(0,\mu]$, and by \eqref{eq:p1_gap}, so does $p^*_1(N)$. By \eqref{eq:zdef}, it follows that $z^*(N)$ converges to
\[
    \bar z:=\frac12\Bigl(\ell+\frac{\mu^2+\sigma^2}{\ell}\Bigr).
\]

It remains to show that $\ell=\pi_\sigma$. To this end, we first show that the sum in \eqref{eq:rev_var_N} converges to $\ln(\bar z/\ell)$. Since $p^*_i\bigl((r^*_i)^2-1\bigr)=2(p^*_{i+1}-p^*_i)+p^*_i(r^*_i-1)^2$, summing over $i$ gives
\[
    \sum_{i=1}^Np^*_i\bigl((r^*_i)^2-1\bigr)=2(z^*-p^*_1)+\sum_{i=1}^Np^*_i(r^*_i-1)^2.
\]
Moreover, combining \eqref{eq:rev_var_N} with \eqref{eq:moment} gives $R^*\sum_{i=1}^Np^*_i\bigl((r^*_i)^2-1\bigr)=\mu^2+\sigma^2-(p^*_1)^2$, and multiplying \eqref{eq:zdef} by $2p^*_1$ gives $\mu^2+\sigma^2-(p^*_1)^2=2p^*_1(z^*-p^*_1)$. Hence,
\begin{align}\label{eq:sq_identity}
    R^*\sum_{i=1}^Np^*_i(r^*_i-1)^2=2(z^*-p^*_1)(p^*_1-R^*).
\end{align}
As $N\to\infty$, the right-hand side of \eqref{eq:sq_identity} converges to zero, while $p^*_i\ge p^*_1$ for all $i$ and $R^*p^*_1\to\ell^2>0$, so that $\sum_{i=1}^N(r^*_i-1)^2\to0$. Since $u-u^2/2\le\ln(1+u)\le u$ for $u\ge0$, and $\sum_{i=1}^N\ln r^*_i=\ln(z^*/p^*_1)$, we have
\[
    0\ \le\ \sum_{i=1}^N(r^*_i-1)-\ln\Bigl(\frac{z^*}{p^*_1}\Bigr)\ \le\ \frac12\sum_{i=1}^N(r^*_i-1)^2,
\]
so that $\sum_{i=1}^N(r^*_i-1)\to\ln(\bar z/\ell)$.

Taking limits in \eqref{eq:rev_var_N}, written as $\mu-p^*_1=R^*\sum_{i=1}^N(r^*_i-1)$, now gives $\mu-\ell=\ell\ln(\bar z/\ell)$, i.e., $\bar z=\ell e^{\mu/\ell-1}$. Inserting this into the definition of $\bar z$, written as $\mu^2+\sigma^2=2\bar z\ell-\ell^2$, gives
\[
    \mu^2+\sigma^2=\ell^2\bigl(2e^{\mu/\ell-1}-1\bigr),
\]
so that $\ell$ solves \eqref{eq:pisigma}. Moreover, $\ell<\mu$, because $2\ell(\bar z-\ell)=\mu^2+\sigma^2-\ell^2>0$ implies $\ln(\bar z/\ell)>0$. Since $\pi_\sigma$ is the unique solution of \eqref{eq:pisigma} in $(0,\mu)$, it follows that $\ell=\pi_\sigma$, and hence $\bar z=\pi_\sigma e^{\mu/\pi_\sigma-1}=z_\infty$.

Finally, we prove that $\bP^*_{\mathrm{seller}}(N) \xrightarrow{d} \mathbb{H}^*_{\mathrm{seller}}$. Let $X_N$ be the corresponding random variable of $\bP^*_{\mathrm{seller}}(N)$ and recall that $\bP^*_{\mathrm{seller}}(N)(X_N=p^*_j(N))=w^*_j$, with $w^*_j$ as in \eqref{eq:var_menu}. Notice that $X_N$ has CDF
\begin{align}
    F_{N}(p) = \begin{cases}
        0, &\quad p < p^*_1,\\
        \sum_{i=1}^jw^*_i, &\quad p^*_j \leq p < p^*_{j+1},\ j=1,\dots,N-1,\\
        1, &\quad p \geq p^*_N.
    \end{cases}
\end{align}
Since $(r^*_i-1)(2z^*-p^*_i-p^*_{i+1})=2z^*(r^*_i-1)-p^*_i\bigl((r^*_i)^2-1\bigr)$ and $p^*_i\bigl((r^*_i)^2-1\bigr)=2(p^*_{i+1}-p^*_i)+p^*_i(r^*_i-1)^2$, applying these identities to the sums in \eqref{eq:var_menu} gives, for $j=1,\dots,N$,
\begin{align}\label{eq:cum_weights}
    \sum_{i=1}^jw^*_i=\frac{2z^*\sum_{i=1}^j(r^*_i-1)-2(p^*_{j+1}-p^*_1)-\sum_{i=1}^jp^*_i(r^*_i-1)^2}{2z^*\sum_{i=1}^N(r^*_i-1)-2(z^*-p^*_1)-\sum_{i=1}^Np^*_i(r^*_i-1)^2}.
\end{align}
Recall that $p^*_1\to\pi_\sigma$, $z^*\to z_\infty$ and $\sum_{i=1}^N(r^*_i-1)^2\to0$ by \eqref{eq:sq_identity}, and that $u-u^2/2\le\ln(1+u)\le u$ for $u\ge0$. Since $\sum_{i=1}^j\ln r^*_i=\ln(p^*_{j+1}/p^*_1)$, it follows that
\[
    0\ \le\ \sum_{i=1}^j\bigl(r^*_i-1\bigr)-\ln\Bigl(\frac{p^*_{j+1}}{p^*_1}\Bigr)\ \le\ \frac12\sum_{i=1}^N\bigl(r^*_i-1\bigr)^2
    \qquad\text{and}\qquad
    0\ \le\ \sum_{i=1}^jp^*_i\bigl(r^*_i-1\bigr)^2\ \le\ z^*\sum_{i=1}^N\bigl(r^*_i-1\bigr)^2
\]
for every $j=1,\dots,N$, so that both middle terms converge to zero, uniformly in $j$. Moreover, with $p^*_{N+1}=z^*$, we have $p^*_{i+1}-p^*_i=p^*_i(r^*_i-1)\le z^*(r^*_i-1)$ for $i=1,\dots,N$, so that
\begin{align}\label{eq:maxgap}
    \max_{1\le i\le N}\bigl(p^*_{i+1}-p^*_i\bigr)\le z^*\Bigl(\sum_{i=1}^N(r^*_i-1)^2\Bigr)^{1/2}\longrightarrow0.
\end{align}
Now fix $p\in(\pi_\sigma,z_\infty)$. By \eqref{eq:maxgap}, $z^*-p^*_N\to0$, so that $p^*_N\to z_\infty$, and hence $p^*_1<p<p^*_N$ for all sufficiently large $N$. Then $F_{N}(p)=\sum_{i=1}^jw^*_i$ for the $j\in\{1,\dots,N-1\}$ with $p^*_j\le p<p^*_{j+1}$, where $p^*_{j+1}\to p$ by \eqref{eq:maxgap}. Together with \eqref{eq:cum_weights}, this gives
\[
    F_{N}(p)\longrightarrow\frac{z_\infty\ln(p/\pi_\sigma)-(p-\pi_\sigma)}{z_\infty\ln(z_\infty/\pi_\sigma)-(z_\infty-\pi_\sigma)}=H^*(p),
\]
where $z_\infty\ln(z_\infty/\pi_\sigma)-(z_\infty-\pi_\sigma)>0$, because $z_\infty>\pi_\sigma$ and $\ln x>1-1/x$ for $x>1$. For $p<\pi_\sigma$, we have $p<p^*_1$ for all sufficiently large $N$, so that $F_{N}(p)=0=H^*(p)$. Similarly, $F_{N}(p)=1=H^*(p)$ for $p>z_\infty$, since then $p>z^*>p^*_N$. Finally, for every $\varepsilon\in(0,z_\infty-\pi_\sigma)$, monotonicity of $F_{N}$ gives $0\le F_{N}(\pi_\sigma)\le F_{N}(\pi_\sigma+\varepsilon)\to H^*(\pi_\sigma+\varepsilon)$, and letting $\varepsilon\downarrow0$ shows that $F_{N}(\pi_\sigma)\to0=H^*(\pi_\sigma)$. In the same way, $F_{N}(z_\infty)\to1=H^*(z_\infty)$. Hence, in the limit when $N\to\infty$, we have $F_{N}(p)\to H^*(p)$ for every $p$, so that $\bP^*_{\mathrm{seller}}(N) \xrightarrow{d} \mathbb{H}^*_{\mathrm{seller}}$.
\end{proof}

\subsection{Mean--MAD convergence}
The next proposition establishes this convergence for the mean--MAD ambiguity set $\cP_d(\mu,d)$, for which the fully randomized revenue $\pi_d$ in \eqref{eq:pid} follows from \citet{Elmachtoub2021value} and \citet{chen2021screening}, together with the convergence of the selling mechanism \eqref{eq:seller_dist}. Here, $W_{-1}$ denotes the lower branch of the Lambert $W$ function, i.e., the inverse of $w\mapsto we^w$ on $(-\infty,-1]$, which is defined on $[-1/e,0)$.
\begin{proposition}\label{prop:asympt}
    Let $\cP=\cP_d(\mu,d)$. The following asymptotic properties hold:
    \begin{enumerate}[label=\textup{(\roman*)}]
        \item\label{it:p1limit} $\lim_{N\to\infty}p^*_1(N) = -\mu\left(1-\frac{d}{2\mu}\right) / W_{-1}\!\left(-\frac{1-\frac{d}{2\mu}}{e}\right) = \pi_d$;
        \item\label{it:Rlimit} $\lim_{N\to\infty}R^*(N) = \pi_d$;
        \item\label{it:distconv} $\bP^*_{\mathrm{seller}}(N) \xrightarrow{d} \mathbb{H}^*_{\mathrm{seller}}$, where $\mathbb{H}^*_{\mathrm{seller}}$ has CDF
        \begin{align}
            H^*(p) =
            \begin{cases}
                0, & p<\pi_d, \\[1mm]
                \dfrac{\ln(p/\pi_d)}{\ln(\mu/\pi_d)}, & p\in[\pi_d,\mu], \\[2mm]
                1, & p>\mu.
            \end{cases}
        \end{align}
    \end{enumerate}
\end{proposition}
\begin{proof}
    First, we prove that $\lim_{N\to\infty}p^*_1(N) = \pi_d$. For $p>0$, define
    $$F_N(p) := \frac{\mu-d/2}{p}-1+N\left[\left(\frac{p}{\mu}\right)^{1/N}-1\right].$$
    For $p \in (0,\mu-d/2]$, notice that
    \begin{align}\label{eq:eq2}
        \lim_{N\to\infty}F_N(p) = \frac{\mu-d/2}{p}-1+\lim_{N\to\infty}\frac{\left[\left(\frac{p}{\mu}\right)^\frac{1}{N}-1\right]}{1/N} = \frac{\mu-d/2}{p}-1 + \ln\left(\frac{p}{\mu}\right) =: F(p)
    \end{align}
    due to L'Hôpital's rule. Since $\lim_{p\downarrow0}F(p) = \infty$, $F(\mu-d/2) = \ln((\mu-d/2)/\mu) < 0$, and $F'(p) = (p-(\mu-d/2))/p^2 < 0$ on the interval $(0,\mu-d/2)$, it follows that $F(p)$ has a unique root on $(0,\mu-d/2]$. 
    
    Next, we show that $\pi_d$ is this unique root. Notice that $F(p)=0$ is equivalent to 
    \begin{align*}
        \ln(p/\mu) = 1-(\mu-d/2)/p &\iff p=\mu e^{1-(\mu-d/2)/p}\\
        &\iff (-(\mu-d/2)/p)e^{-(\mu-d/2)/p} = -(\mu-d/2)/(\mu e)\\
        &\iff (-(\mu-d/2)/p)=W_{-1}\left(-\frac{\mu-d/2}{\mu e}\right)\\
        &\iff p = -\left(\mu-\frac{d}{2}\right) / W_{-1}\!\left(-\frac{1-\frac{d}{2\mu}}{e}\right) = \pi_d.
    \end{align*}
    Since $-\frac{1-\frac{d}{2\mu}}{e} \in [-\frac{1}{e},0)$, both branches $W_{-1}$ and $W_{0}$ are defined. However, since $W_0$ would give $p > \mu-d/2$ while $W_{-1}$ gives $p \leq \mu-d/2$, we select $W_{-1}$.

    We will now show that $p^*_1(N) \to \pi_d$ when $N \to\infty$. Let $\varepsilon > 0$ sufficiently small. Then $F(\pi_d-\varepsilon) > 0$ and $F(\pi_d+\varepsilon) < 0$, and because $\lim_{N\to\infty}F_N(p) = F(p)$ for fixed $p$, we obtain $\lim_{N\to\infty}F_N(\pi_d-\varepsilon)=F(\pi_d-\varepsilon)>0$ and $\lim_{N\to\infty}F_N(\pi_d+\varepsilon)=F(\pi_d+\varepsilon)<0$. Hence, there exists a $N_0 \in \N$ such that for all $N > N_0$ it holds that
    $F_N(\pi_d-\varepsilon) > 0$ and $F_N(\pi_d+\varepsilon) < 0$. Since $F_N(p)$ is continuous, the intermediate value theorem states that a root lies in $(\pi_d-\varepsilon,\pi_d+\varepsilon)$, and since $p^*_1(N)$ is the unique root in $(0,\mu-d/2]$, we know that $p^*_1(N) \in (\pi_d-\varepsilon,\pi_d+\varepsilon)$ for $N > N_0$. Therefore, $|p^*_1(N) - \pi_d| < \varepsilon$, and since $\varepsilon$ was chosen arbitrarily, this is equivalent to $\lim_{N\to\infty}p^*_1(N)=\pi_d$.

    Second, we prove that $\lim_{N\to\infty}R^*(N)=\pi_d$. Notice that \begin{align}
        \lim_{N\to\infty}R^*(N) = \lim_{N\to\infty}\frac{N+1}{N}\lim_{N\to\infty}p^*_1(N)-\lim_{N\to\infty}\frac{\mu-d/2}{N}=\pi_d.
    \end{align}

    Finally, we prove that $\bP^*_{\mathrm{seller}}(N) \xrightarrow{d} \mathbb{H}^*_{\mathrm{seller}}$. Let $X_N$ be the corresponding random variable of $\bP^*_{\mathrm{seller}}(N)$ and recall that $\bP^*_{\mathrm{seller}}(N)(X_N=p^*_j(N))=1/N$. Now define $Y_N = \ln(X_N)$ and notice $Y_N$ takes values
    $$\ln(p^*_j(N))=\ln(p^*_1(N))+\frac{j-1}{N}\ln\left(\frac{\mu}{p^*_1(N)}\right)$$ with probability $1/N$ for $j = 1,\dots,N$. Alternatively, we can write $$Y_N = \ln(p^*_1(N))+U_N\ln\left(\frac{\mu}{p^*_1(N)}\right),$$
    where $U_N$ is (discrete) uniformly distributed on $\{0,\frac{1}{N},\dots,\frac{N-1}{N}\}$. 
    
    We now show that $U_N \xrightarrow{d} U\sim \mathrm{Unif}(0,1)$. Notice that $U_N$ has CDF
    \begin{align}
        F_{N,U}(u) = \begin{cases}
            0, &\quad u < 0,\\
            \frac{\lfloor Nu\rfloor + 1}{N}, &\quad 0 \leq u < \frac{N-1}{N},\\
            1, &\quad u \geq \frac{N-1}{N}.
        \end{cases}
    \end{align}
    It follows that for $u \in [0,1)$ we have $u < F_{N,U}(u) \leq u +\frac{1}{N}$, so that in the limit when $N \to \infty$ we have $F_{N,U}(u) \to F_U(u)$ with $F_U(u)$ the CDF of the uniform distribution on $(0,1)$.

    Next, we apply Slutsky's theorem \citep{bain1992introduction} to conclude that \begin{align*}
        Y_N \xrightarrow{d} \ln(\pi_d) + U\ln\left(\frac{\mu}{\pi_d}\right).
    \end{align*}
    Then, due to the continuous mapping theorem \citep{bain1992introduction} with $g(Y)=e^{Y}$ we obtain
    \begin{align*}
        X_N \xrightarrow{d} e^{\ln(\pi_d) + U\ln\left(\frac{\mu}{\pi_d}\right)} \sim \mathbb{H}^*_{\mathrm{seller}},
    \end{align*}
    since, \begin{align}
         \bP\left(e^{\ln(\pi_d) + U\ln\left(\frac{\mu}{\pi_d}\right)} \leq p\right) =  \bP\left(\ln(\pi_d) + U\ln\left(\frac{\mu}{\pi_d}\right) \leq \ln(p) \right) =  \bP\left(U \leq \dfrac{\ln(p/\pi_d)}{\ln(\mu/\pi_d)} \right).
    \end{align}
    This corresponds with CDF $H^*(p)$, since the boundary points (of the support) $p^*_1(N) \to \pi_d$ and $p^*_N(N) \to \mu$ when $N \to \infty$.
\end{proof}

\section{Proof of Theorem \ref{th:convergence_rate}}\label{ap:convergence_rate}
We prove the two parts of Theorem~\ref{th:convergence_rate} separately, following the same two steps. We first derive a $1/N$-expansion of $R^*(N)$ from the optimal menus of Theorems~\ref{th:var_N} and~\ref{th:mad_N}. We then show that the leading coefficient is increasing in the dispersion parameter.

\subsection{Mean--variance asymptotics}\label{ap:convergence_rate_var}
Recall from the proof of Proposition~\ref{prop:asympt_var} that $p^*_1\to\pi_\sigma$, $R^*\to\pi_\sigma$ and $z^*\to z_\infty$, where
\[
    z_\infty=\frac12\Bigl(\pi_\sigma+\frac{\mu^2+\sigma^2}{\pi_\sigma}\Bigr)=\pi_\sigma e^{\mu/\pi_\sigma-1},
\]
so that $\mu-\pi_\sigma=\pi_\sigma\ln(z_\infty/\pi_\sigma)$, $\mu^2+\sigma^2=2\pi_\sigma z_\infty-\pi_\sigma^2$ and $z_\infty>\pi_\sigma$. By \eqref{eq:p1_gap} and \eqref{eq:sq_identity}, we have
\[
    p^*_1-R^*=O\Bigl(\frac1N\Bigr),\qquad\sum_{i=1}^N(r^*_i-1)^2=O\Bigl(\frac1N\Bigr).
\]
Recall from the proof of Proposition~\ref{prop:asympt_var} that $0\le\sum_{i=1}^N(r^*_i-1)-\ln(z^*/p^*_1)\le\frac12\sum_{i=1}^N(r^*_i-1)^2$, so that $\sum_{i=1}^N(r^*_i-1)=\ln(z^*/p^*_1)+O(1/N)$. Substituting this in \eqref{eq:rev_var_N}, written as $\mu-p^*_1=R^*\sum_{i=1}^N(r^*_i-1)$, and using $R^*=p^*_1+O(1/N)$ results in
\[
    p^*_1\Bigl(1+\ln\frac{z^*}{p^*_1}\Bigr)-\mu=O\Bigl(\frac1N\Bigr).
\]
Then, substituting $z^*=\frac12\bigl(p^*_1+(\mu^2+\sigma^2)/p^*_1\bigr)$ from \eqref{eq:zdef} and using the Taylor expansion around $p^*_1=\pi_\sigma$, we obtain
\[
    \Bigl(\ln\frac{z_\infty}{\pi_\sigma}-1+\frac{\pi_\sigma}{z_\infty}+o(1)\Bigr)(p^*_1-\pi_\sigma)=O\Bigl(\frac1N\Bigr).
\]
Since $\ln x>1-1/x$ for $x>1$, the coefficient of $p^*_1-\pi_\sigma$ is bounded away from zero for all sufficiently large $N$. Hence, by combining this with \eqref{eq:zdef} and \eqref{eq:p1_gap}, we conclude that
\begin{align}\label{eq:endpoints}
    p^*_1=\pi_\sigma+O\Bigl(\frac1N\Bigr),\qquad z^*=z_\infty+O\Bigl(\frac1N\Bigr),\qquad R^*=\pi_\sigma+O\Bigl(\frac1N\Bigr).
\end{align}
In particular, all prices lie in $[\pi_\sigma/2,2z_\infty]$ for all sufficiently large $N$.

Next, recall from \eqref{eq:Psi} that $R^*=A(p^*,z^*)/D(p^*,z^*)$, where, with $p_{N+1}=z$ and $r_i=p_{i+1}/p_i$,
\[
    A(p,z)=(z-p_1)^2-(z-\mu)^2-\sigma^2,\qquad D(p,z)=\sum_{i=1}^N(r_i-1)(2z-p_i-p_{i+1})=\sum_{i=1}^N\frac{(z-p_i)^2-(z-p_{i+1})^2}{p_i}.
\]
To expand $D(p^*,z^*)$, define, for $0<s\le z$,
\[
    J(s,z):=\sqrt z\ln\frac{\sqrt z+\sqrt{z-s}}{\sqrt z-\sqrt{z-s}}-2\sqrt{z-s},
\]
so that $J(z,z)=0$ and $\partial J(s,z)/\partial s=-\sqrt{z-s}/s$. We will now show that
\begin{align}\label{eq:D_exp}
    D(p^*,z^*)=2z^*\ln\frac{z^*}{p^*_1}-2(z^*-p^*_1)+\frac{J(p^*_1,z^*)^2}{N}+O\Bigl(\frac{\ln N}{N^2}\Bigr).
\end{align}
Consider any points $p^*_1=q_1<q_2<\dots<q_{N+1}=z^*$, let $\delta_i:=q_{i+1}-q_i$, and define
\[
    J_i:=J(q_i,z^*)-J(q_{i+1},z^*),\qquad e_i:=2z^*\Bigl(\frac{\delta_i}{q_i}-\ln\Bigl(1+\frac{\delta_i}{q_i}\Bigr)\Bigr)-\frac{\delta_i^2}{q_i}.
\]
Recall from \eqref{eq:Psi} that the $i$-th term of $D(q,z^*)$ equals $\bigl((z^*-q_i)^2-(z^*-q_{i+1})^2\bigr)/q_i=2z^*\delta_i/q_i-2\delta_i-\delta_i^2/q_i$, which in turn equals $e_i+2z^*\ln(q_{i+1}/q_i)-2\delta_i$. Summing over $i$ gives
\[
    D(q,z^*)=2z^*\ln\frac{z^*}{p^*_1}-2(z^*-p^*_1)+\sum_{i=1}^Ne_i,\qquad\sum_{i=1}^NJ_i=J(p^*_1,z^*).
\]
Using the Taylor expansion $x-\ln(1+x)=x^2/2+O(x^3)$ around $x=0$ with $x=\delta_i/q_i$, we obtain $e_i=(z^*-q_i)\delta_i^2/q_i^2+O(\delta_i^3)$. Furthermore, by the mean value theorem, there exists a $\xi_i\in(q_i,q_{i+1})$ such that $J_i=\delta_i\sqrt{z^*-\xi_i}/\xi_i$, and since $t\mapsto(z^*-t)/t^2$ has a bounded derivative on $[\pi_\sigma/2,2z_\infty]$, also $J_i^2=(z^*-q_i)\delta_i^2/q_i^2+O(\delta_i^3)$. Hence, $|e_i-J_i^2|\le C\delta_i^3$ for a constant $C>0$ that does not depend on $N$ or $i$.

We also need lower bounds on $e_i$ and $J_i$. As functions of $\delta_i\in[0,z^*-q_i]$, both equal zero at $\delta_i=0$, and since $q_i(q_i+\delta_i)\le(z^*)^2$, their derivatives with respect to $\delta_i$ satisfy
\[
    \frac{2\delta_i(z^*-q_i-\delta_i)}{q_i(q_i+\delta_i)}\ge\frac{2\delta_i(z^*-q_i-\delta_i)}{(z^*)^2},\qquad\frac{\sqrt{z^*-q_i-\delta_i}}{q_i+\delta_i}\ge\frac{\sqrt{z^*-q_i-\delta_i}}{z^*}.
\]
The right-hand sides are the derivatives with respect to $\delta_i$ of $\frac2{(z^*)^2}\bigl(\frac12(z^*-q_i)\delta_i^2-\frac13\delta_i^3\bigr)$ and $\frac2{3z^*}\bigl((z^*-q_i)^{3/2}-(z^*-q_i-\delta_i)^{3/2}\bigr)$, which also equal zero at $\delta_i=0$. Hence, using $\delta_i\le z^*-q_i$ and $(z^*-q_i-\delta_i)^{3/2}\le(z^*-q_i-\delta_i)\sqrt{z^*-q_i}$, we obtain
\[
    e_i\ge\frac{(z^*-q_i)\delta_i^2}{3(z^*)^2},\qquad J_i\ge\frac{2\delta_i\sqrt{z^*-q_i}}{3z^*}.
\]
Consequently, after increasing $C$ if necessary,
\begin{align}\label{eq:ratio}
    \max\Bigl\{\frac{e_i}{J_i^2},\frac{J_i^2}{e_i}\Bigr\}\le1+\frac{C\delta_i}{z^*-q_i}.
\end{align}
If the $\delta_i$ are nondecreasing in $i$, then $\delta_N\ge(z^*-p^*_1)/N$. Since $1-1/x\le\ln x$ for $x>0$, we have $\delta_i/(z^*-q_i)\le\ln\bigl((z^*-q_i)/(z^*-q_{i+1})\bigr)$ for $i<N$, which by summing results in
\begin{align}\label{eq:logsum}
    \sum_{i=1}^N\frac{\delta_i}{z^*-q_i}\le1+\ln\frac{z^*-p^*_1}{\delta_N}\le1+\ln N.
\end{align}

The gaps $p^*_{i+1}-p^*_i$ are increasing in $i$, since $p^*_{i+1}-p^*_i=p^*_i(r^*_i-1)$. Then, by the Cauchy--Schwarz inequality, \eqref{eq:ratio} and \eqref{eq:logsum}, we obtain
\[
    J(p^*_1,z^*)^2=\Bigl(\sum_{i=1}^NJ_i\Bigr)^2\le\Bigl(\sum_{i=1}^Ne_i\Bigr)\Bigl(\sum_{i=1}^N\frac{J_i^2}{e_i}\Bigr)\le\Bigl(\sum_{i=1}^Ne_i\Bigr)\bigl(N+C(1+\ln N)\bigr).
\]
For the reverse bound, choose the points $q_i$ such that $J_i=J(p^*_1,z^*)/N$ for every $i$. Such points exist since $J(\cdot,z^*)$ is continuous and strictly decreasing on $[p^*_1,z^*]$ with $J(z^*,z^*)=0$. Moreover, their gaps are nondecreasing, since $J_i=\delta_i\sqrt{z^*-\xi_i}/\xi_i$ with $\xi_1<\xi_2<\dots<\xi_N$ and $\sqrt{z^*-t}/t$ is decreasing in $t$. Hence, \eqref{eq:ratio} and \eqref{eq:logsum} imply
\[
    \sum_{i=1}^Ne_i\le\frac{J(p^*_1,z^*)^2}{N^2}\bigl(N+C(1+\ln N)\bigr).
\]
Recall from the proof of Theorem~\ref{th:var_N} that $(p^*,z^*)$ maximizes $\Psi(p,z)=A(p,z)/D(p,z)$ with $A(p^*,z^*)>0$, where $A(p,z)$ depends on $p$ only through $p_1$. Hence, the optimal prices minimize $D(\cdot,z^*)$ among all prices with lowest price $p^*_1$, so that the sum of the $e_i$ for the optimal prices lies between
\[
    \frac{J(p^*_1,z^*)^2}{N+C(1+\ln N)}\qquad\text{and}\qquad\frac{J(p^*_1,z^*)^2}{N^2}\bigl(N+C(1+\ln N)\bigr).
\]
Both equal $J(p^*_1,z^*)^2/N+O(\ln N/N^2)$, which proves \eqref{eq:D_exp}.

We now extract the expansion of $R^*$. Recall that $A(p^*,z^*)=2z^*(\mu-p^*_1)-\bigl(\mu^2+\sigma^2-(p^*_1)^2\bigr)$. By the identities for $\mu-\pi_\sigma$ and $\mu^2+\sigma^2$ above, the Taylor expansions of $A(p^*,z^*)$ and of $\pi_\sigma\bigl(2z^*\ln(z^*/p^*_1)-2(z^*-p^*_1)\bigr)$ around $(p^*_1,z^*)=(\pi_\sigma,z_\infty)$ have the same constant and first-order terms. Hence, by \eqref{eq:endpoints},
\[
    A(p^*,z^*)=\pi_\sigma\Bigl(2z^*\ln\frac{z^*}{p^*_1}-2(z^*-p^*_1)\Bigr)+O\Bigl(\frac1{N^2}\Bigr).
\]
Since $R^*=\Psi(p^*,z^*)$, we have $R^*D(p^*,z^*)=A(p^*,z^*)$, and subtracting $\pi_\sigma$ times \eqref{eq:D_exp} results in
\[
    (R^*-\pi_\sigma)D(p^*,z^*)=-\pi_\sigma\frac{J(p^*_1,z^*)^2}{N}+O\Bigl(\frac{\ln N}{N^2}\Bigr).
\]
By \eqref{eq:endpoints} and \eqref{eq:D_exp}, $D(p^*,z^*)$ and $J(p^*_1,z^*)^2$ differ from $2z_\infty\bigl(\ln(z_\infty/\pi_\sigma)-1+\pi_\sigma/z_\infty\bigr)>0$ and $J(\pi_\sigma,z_\infty)^2$ by $O(1/N)$. Therefore,
\[
    R^*(N)=\pi_\sigma\Bigl(1-\frac{\kappa_\sigma}N+O\Bigl(\frac{\ln N}{N^2}\Bigr)\Bigr),\qquad\kappa_\sigma=\frac{J(\pi_\sigma,z_\infty)^2}{2z_\infty\bigl(\ln(z_\infty/\pi_\sigma)-1+\pi_\sigma/z_\infty\bigr)}.
\]
Since $a_\sigma=1-\pi_\sigma/z_\infty$, we have $z_\infty-\pi_\sigma=a_\sigma z_\infty$, so that the definition of $J$ gives
\[
    J(\pi_\sigma,z_\infty)=2\sqrt{z_\infty}\Bigl(\frac12\ln\frac{1+\sqrt{a_\sigma}}{1-\sqrt{a_\sigma}}-\sqrt{a_\sigma}\Bigr),
\]
while $\ln(z_\infty/\pi_\sigma)-1+\pi_\sigma/z_\infty=\mu/\pi_\sigma-1-a_\sigma$. This gives the stated expression for $\kappa_\sigma$.

It remains to show that $\kappa_\sigma$ is increasing in $\sigma$. For $a\in(0,1)$, define
\[
    H(a):=\frac12\ln\frac{1+\sqrt a}{1-\sqrt a}-\sqrt a,\qquad L(a):=-\ln(1-a)-a.
\]
Since $\mu/\pi_\sigma-1=-\ln(1-a_\sigma)$, we have $\kappa_\sigma=2H(a_\sigma)^2/L(a_\sigma)$. Clearly, $H(a)$ and $L(a)$ converge to zero as $a\downarrow0$, and since $H'(a)=\sqrt a/(2(1-a))>0$ and $L'(a)=a/(1-a)>0$, both are positive on $(0,1)$. Differentiating gives
\[
    \frac{d}{da}\,\kappa_\sigma=\frac{2\sqrt a\,H(a)}{(1-a)L(a)^2}\bigl(L(a)-\sqrt a\,H(a)\bigr),
\]
where, using the power series of $-\ln(1-a)$ and of $\ln\bigl((1+\sqrt a)/(1-\sqrt a)\bigr)$,
\[
    L(a)-\sqrt a\,H(a)=-\ln(1-a)-\frac{\sqrt a}2\ln\frac{1+\sqrt a}{1-\sqrt a}=\sum_{k=2}^\infty\frac{k-1}{k(2k-1)}\,a^k>0.
\]
Hence, $\kappa_\sigma$ is increasing in $a_\sigma$. Moreover, $a_\sigma=1-e^{1-\mu/\pi_\sigma}$ is decreasing in $\pi_\sigma$, so it remains to show that $\pi_\sigma$ is decreasing in $\sigma$. By Proposition~\ref{prop:prices_var}, $R^*(N)$ is decreasing in $\sigma$ for every $N$, so that its limit $\pi_\sigma$ is nonincreasing in $\sigma$. The right-hand side of~\eqref{eq:pisigma} depends on $\sigma$ only through $\pi_\sigma$, while the left-hand side is strictly increasing in $\sigma$. Hence, $\pi_\sigma$ cannot take the same value for two different values of $\sigma$, and is therefore decreasing in $\sigma$. We conclude that $\kappa_\sigma$ is increasing in $\sigma$.

\subsection{Mean--MAD asymptotics}\label{ap:convergence_rate_mad}
Define $t_N = \ln\left(\mu/p_1(N)\right)$, so that $p_1(N) = \mu e^{-t_N}$. Substituting this in \eqref{eq:p1} results in
    \begin{align}
        \left(1-\frac{d}{2\mu}\right)e^{t_N}=1+N\left(1-e^{-t_N/N}\right).
    \end{align}
    Then, using the Taylor expansion of $e^x$ around $x=0$ with $x=-t_N/N$, we obtain
    \begin{align}
        e^{-t_N/N} = 1 - \frac{t_N}{N} + \frac{t_N^2}{2N^2}+O\left(\frac{1}{N^3}\right).
    \end{align}
    Therefore,
    \begin{align}\label{eq:eq3}
        \left(1-\frac{d}{2\mu}\right)e^{t_N}=1+t_N-\frac{t_N^2}{2N}+O\left(\frac{1}{N^2}\right).
    \end{align}

    Next, using Proposition \ref{prop:asympt} we have that $\lim_{N\xrightarrow{}\infty} t_N = \ln(\mu/\pi_d) =: t_{\infty}$. Substituting this in \eqref{eq:eq2} results in
    \begin{align}\label{eq:eq4}
        \left(1-\frac{d}{2\mu}\right)e^{t_\infty}=1+t_{\infty}.
    \end{align}
    Then, subtracting \eqref{eq:eq4} from \eqref{eq:eq3} we obtain
    \begin{align}\label{eq:eq5}
        \left(1-\frac{d}{2\mu}\right)(e^{t_N}-e^{t_\infty}) = t_N-t_{\infty}-\frac{t_N^2}{2N}+O\left(\frac{1}{N^2}\right).
    \end{align}

    We will now show that $t_N-t_\infty = O(1/N)$. By the mean value theorem, there exists a $\xi_N \in (t_N,t_{\infty})$ such that
    \begin{align}
        e^{t_N}-e^{t_{\infty}}=e^{\xi_N}(t_N - t_\infty),
    \end{align}
    so that \eqref{eq:eq5} becomes
    \begin{align}\label{eq:eq6}
        \left(\left(1-\frac{d}{2\mu}\right)e^{\xi_N}-1\right)(t_N-t_{\infty}) = -\frac{t_N^2}{2N}+O\left(\frac{1}{N^2}\right) = O\left(\frac{1}{N}\right),
    \end{align}
    since $t_N^2 \xrightarrow{}t^2_\infty \in (0,M]$ with $M$ some sufficiently large positive constant. Furthermore, since $t_N \xrightarrow{} t_{\infty}$, also $\xi_N \xrightarrow{} t_\infty$. Therefore, due to \eqref{eq:eq4}, we obtain
    \begin{align}
        \left(\left(1-\frac{d}{2\mu}\right)e^{\xi_N}-1\right) \xrightarrow{} \left(\left(1-\frac{d}{2\mu}\right)e^{t_{\infty}}-1\right) = t_{\infty} \in (0,M].
    \end{align}
    Hence, by combining this with \eqref{eq:eq6}, we conclude that
   \begin{align}
        t_N-t_{\infty} = O\left(\frac{1}{N}\right).
    \end{align}
    In turn, this implies that \eqref{eq:eq6} becomes
    \begin{align}
        &\left(t_{\infty}+O(1/N)\right)(t_N-t_{\infty})=-\frac{t^2_\infty}{2N}+O(1/N^2)\\
        &\iff t_{\infty}(t_N-t_{\infty}) + O(1/N^2) = -\frac{t^2_\infty}{2N}+O(1/N^2)\\
        &\iff t_N-t_{\infty} = -\frac{t_\infty}{2N}+O(1/N^2).
    \end{align}
    
    Next, notice that
    \begin{align}\label{eq:eq7}
        p_1(N) = \mu e^{-t_N}=\mu e^{-t_{\infty}}e^{-(t_N-t_\infty)}=\pi_de^{-(t_N-t_\infty)}.
    \end{align}
    Hence, using the Taylor expansion of $e^{-x}$ around $x=0$ with $x=t_N-t_\infty$, we obtain
    \begin{align}
        p_1(N) = \pi_d\left(1-(t_N-t_{\infty})+O(1/N^2)\right) = \pi_d\left(1+\frac{\ln(\mu/\pi_d)}{2N}+O(1/N^2)\right).
    \end{align}
    Therefore, plugging \eqref{eq:eq7} into \eqref{eq:rev_N_2a}, and recalling \eqref{eq:eq4}, we obtain
    \begin{align}
        R^*(N) &= \frac{N+1}{N}\pi_d\left(1+\frac{\ln(\mu/\pi_d)}{2N}+O(1/N^2)\right)-\frac{\mu-d/2}{N}\\
        &=\left(1+\frac{1}{N}\right)\pi_d\left(1+\frac{\ln(\mu/\pi_d)}{2N}+O(1/N^2)\right)-\frac{\pi_d(1+\ln(\mu/\pi_d))}{N}\\
        &= \pi_d\left(1+\frac{1}{N}+\frac{\ln(\mu/\pi_d)}{2N}+\frac{\ln(\mu/\pi_d)}{2N^2}+O(1/N^2)\right)-\frac{\pi_d(1+\ln(\mu/\pi_d))}{N}\\
        &= \pi_d\left(1+\frac{1}{N}+\frac{\ln(\mu/\pi_d)}{2N}+O(1/N^2)\right)-\frac{\pi_d(1+\ln(\mu/\pi_d))}{N}\\
        &= \pi_d\left(1-\frac{\ln(\mu/\pi_d)}{2N}+O(1/N^2)\right).
    \end{align}

    It remains to show that the leading coefficient $\ln(\mu/\pi_d)/2$ is increasing in $d$, i.e., that $\pi_d$ is decreasing in $d$. By Proposition~\ref{prop:prices_decrease_MAD}, $R^*(N) = p_1^* p_N^*/\mu$ is decreasing in $d$ for every $N$, so that its limit $\pi_d$ is nonincreasing in $d$. The left-hand side of~\eqref{eq:pid} is strictly decreasing in $d$, while the right-hand side depends on $d$ only through $\pi_d$. Hence, $\pi_d$ cannot take the same value for two different values of $d$, and is therefore decreasing in $d$. We conclude that $\kappa_d$ is increasing in $d$.

\section{Proofs of Lemmas \ref{lem:strong_dualityvar},
\ref{lem:strong_duality}, and \ref{lem:strong_duality2}}\label{ap:strong_duality}

\begin{proof}
We prove the three strong-duality results jointly. In each case, write
$\mathcal X$ for the support of the primal measure, $h(x)$ for the vector
of moment functions, and $b$ for the prescribed moment vector, as summarized
in Table \ref{tab:strong_duality_moments}.

\begin{table}[H]
\centering
\begin{tabular}{c|c|c|c}
\hline
Ambiguity set & $\mathcal X$ & $h(x)$ & $b$ \\ \hline
Mean--support
& $[0,\mu]$
& $(1,x)$
& $(1,\mu-d/2)$ \\[1mm]

Mean--MAD
& $\R_+$
& $(1,x,|x-\mu|)$
& $(1,\mu,d)$ \\[1mm]

Mean--variance
& $\R_+$
& $(1,x,x^2)$
& $(1,\mu,\mu^2+\sigma^2)$ \\ \hline
\end{tabular}
\caption{Parameters for the different ambiguity sets.}
\label{tab:strong_duality_moments}
\end{table}

For the problem under consideration, define the value function
$$
    \phi(y)
    :=
    \inf_{\bP\in\mathcal M_+(\mathcal X)}
    \left\{
        \int_{\mathcal X}G(x)\,{\rm d}\bP(x):
        \int_{\mathcal X}h(x)\,{\rm d}\bP(x)=y
    \right\},
$$
where $\mathcal M_+(\mathcal X)$ denotes the non-negative measures on
$\mathcal X$. Thus the corresponding primal value equals $\phi(b)$. We now establish the
properties of $\phi$ needed to construct a dual solution attaining this value. We first show that $b$ is an interior point of the moment cone
and that $\phi$ is finite on a neighborhood of $b$. Together with convexity of $\phi$,
this allows us to apply the supporting-hyperplane theorem to the epigraph of
$\phi$ at $(b,\phi(b))$, giving a supporting vector $\lambda^*$. Positive
homogeneity of $\phi$ will then imply
\[
    (\lambda^*)^\top b=\phi(b),
\]
and evaluating the resulting supporting-hyperplane inequality at the moment vector
$h(x)$ will show that $\lambda^*$ is dual feasible. This will establish
strong duality.

We proceed with verifying that $b$ is an interior point of the moment cone
$$
    \mathcal C
    :=
    \left\{
        \int_{\mathcal X}h(x)\,{\rm d}\bP(x):
        \bP\in\mathcal M_+(\mathcal X)
    \right\}.
$$
For mean--support, $0<d<2\mu$ gives
$$
    (1,\mu-d/2)
    =
    \frac{d}{2\mu}h(0)
    +
    \left(1-\frac{d}{2\mu}\right)h(\mu),
$$
and both coefficients are strictly positive. Since $h(0)$ and $h(\mu)$
are linearly independent, $b$ lies in the interior of $\mathcal C$.

For mean--MAD, choose $y>\mu$ sufficiently large and define
$$
    \alpha_0=\frac{d}{2\mu},
    \qquad
    \alpha_y=\frac{d}{2(y-\mu)},
    \qquad
    \alpha_\mu=1-\alpha_0-\alpha_y.
$$
Because $0<d<2\mu$, $y$ can be chosen so that
$\alpha_0,\alpha_\mu,\alpha_y>0$. A direct calculation gives
$$
    (1,\mu,d)
    =
    \alpha_0h(0)+\alpha_\mu h(\mu)+\alpha_yh(y).
$$
Moreover, $h(0)$, $h(\mu)$, and $h(y)$ are linearly independent, since
their determinant equals $2\mu(y-\mu)>0$. Hence $b$ is again an interior
point of $\mathcal C$.

For mean--variance, choose
$$
    y>\mu+\frac{\sigma^2}{\mu}
$$
and define
$$
    \alpha_0=\frac{\sigma^2}{\mu y},
    \qquad
    \alpha_y=\frac{\sigma^2}{y(y-\mu)},
    \qquad
    \alpha_\mu=1-\alpha_0-\alpha_y.
$$
The choice of $y$ ensures that all three coefficients are strictly
positive, and
$$
    (1,\mu,\mu^2+\sigma^2)
    =
    \alpha_0h(0)+\alpha_\mu h(\mu)+\alpha_yh(y).
$$
The three vectors are linearly independent, with determinant
$\mu y(y-\mu)>0$, so $b$ is an interior point of $\mathcal C$.

Furthermore, $\phi$ is positively homogeneous. Indeed, $\phi(0)=0$, and by linearity of integration with respect to the measure, scaling $P$ by $t>0$ scales both its moment vector and objective value by $t$, so
\[
    \phi(ty)=t\phi(y), \qquad t>0.
\]

In addition, $\phi$ is convex. To see this, let
$y_1,y_2\in\mathcal C$ and $\theta\in[0,1]$. For any $\varepsilon>0$, choose
$P_1,P_2\in\mathcal M_+(\mathcal X)$ feasible for $y_1$ and $y_2$, respectively,
such that
\[
    \int G\,dP_i \leq \phi(y_i)+\varepsilon,
    \qquad i=1,2.
\]
Then $P_\theta:=\theta P_1+(1-\theta)P_2$ is feasible for
$\theta y_1+(1-\theta)y_2$, since
\[
    \int h\,dP_\theta
    =
    \theta\int h\,dP_1+(1-\theta)\int h\,dP_2
    =
    \theta y_1+(1-\theta)y_2.
\]
Hence
\[
\begin{aligned}
    \phi\bigl(\theta y_1+(1-\theta)y_2\bigr)
    &\leq \int G\,dP_\theta \\
    &= \theta\int G\,dP_1+(1-\theta)\int G\,dP_2 \\
    &\leq \theta\phi(y_1)+(1-\theta)\phi(y_2)+\varepsilon.
\end{aligned}
\]
Letting $\varepsilon\downarrow0$ proves convexity.

It remains to verify that
$\phi$ is finite on an open neighborhood of $b$. For each of the three cases, let $x_1,\ldots,x_m$ denote the support points
used in the preceding representation of $b$, and define
\[
    K:=\operatorname{cone}\{h(x_1),\ldots,h(x_m)\}.
\]
The same argument used above for $\mathcal C$ shows that $b$ is an interior
point of $K$. Hence there exists
$\varepsilon>0$ such that
\[
    B_\varepsilon(b)\subseteq K.
\]
Thus, for every $y\in B_\varepsilon(b)$, there exist coefficients
$\gamma_1,\ldots,\gamma_m\geq0$ such that
\[
    y=\sum_{i=1}^m \gamma_i h(x_i).
\]
Consequently, the finite-support measure
\[
    P_y:=\sum_{i=1}^m \gamma_i\delta_{x_i}
\]
satisfies
\[
    \int h\,dP_y
    =
    \sum_{i=1}^m \gamma_i h(x_i)
    =
    y,
\]
and is therefore feasible for the problem defining $\phi(y)$. It follows that
\[
    \phi(y)
    \leq
    \int G\,dP_y
    =
    \sum_{i=1}^m \gamma_i G(x_i)
    <\infty.
\]
Hence $\phi$ is finite from above on an open neighborhood of $b$. Moreover, $\phi(y)>-\infty$: in the mean--support and
mean--MAD cases, $G(x)\geq G(0)>-\infty$, while in the mean--variance case,
$G(x)\geq0$.

Although $b$ is an interior point of the moment cone $\mathcal C$, the point
$(b,\phi(b))$ lies on the boundary of the epigraph of $\phi$, since it lies
on the graph of $\phi$. Since $\phi$ is convex, its epigraph is convex.
Moreover, because $\phi$ is finite on an open neighborhood of $b$, the
supporting hyperplane at $(b,\phi(b))$ may be represented as the graph of an
affine function. Hence, by the supporting-hyperplane theorem, there exists a
vector $\lambda^*$ such that
\[
    \phi(y)
    \geq
    \phi(b)+(\lambda^*)^\top(y-b),
    \qquad y\in\mathcal C.
\]
Positive homogeneity gives $\phi(0)=0$ and $\phi(2b)=2\phi(b)$. Applying
the preceding inequality first at $y=0$ and then at $y=2b$ yields
$$
    (\lambda^*)^\top b=\phi(b).
$$
Finally, for every $x\in\mathcal X$, the point mass $\delta_x$ is
feasible for the moment vector $h(x)$, and therefore
$$
    G(x)
    \geq
    \phi(h(x))
    \geq
    \phi(b)+(\lambda^*)^\top(h(x)-b)
    =
    (\lambda^*)^\top h(x).
$$
Thus $\lambda^*$ is feasible for the corresponding dual problem and its
dual objective equals the primal value $\phi(b)$. Weak duality gives the
reverse inequality, so the primal and dual values coincide and the dual
optimum is attained. Applying this argument to the three rows of the table proves Lemmas
\ref{lem:strong_duality}, \ref{lem:strong_duality2}, and
\ref{lem:strong_dualityvar}.
\end{proof}

\section{Decreasing optimal prices under mean--MAD}\label{ap:MADdecrease}
\begin{figure}[h!]
  \centering
  \begin{subfigure}[b]{0.48\linewidth}
    \centering
    \caption*{$N = 1$ (deterministic pricing)}
    \includegraphics[width=\linewidth]{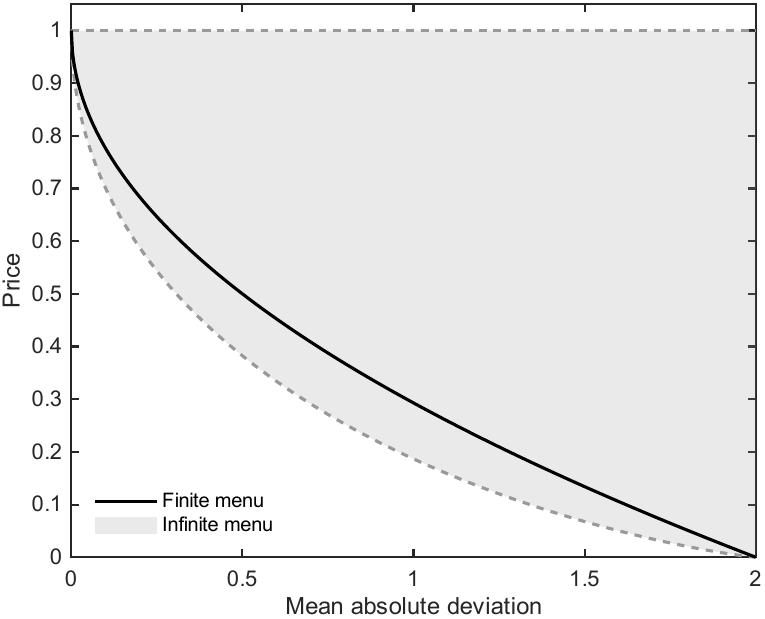}
  \end{subfigure}\hfill
  \begin{subfigure}[b]{0.48\linewidth}
    \centering
    \caption*{$N = 2$}
    \includegraphics[width=\linewidth]{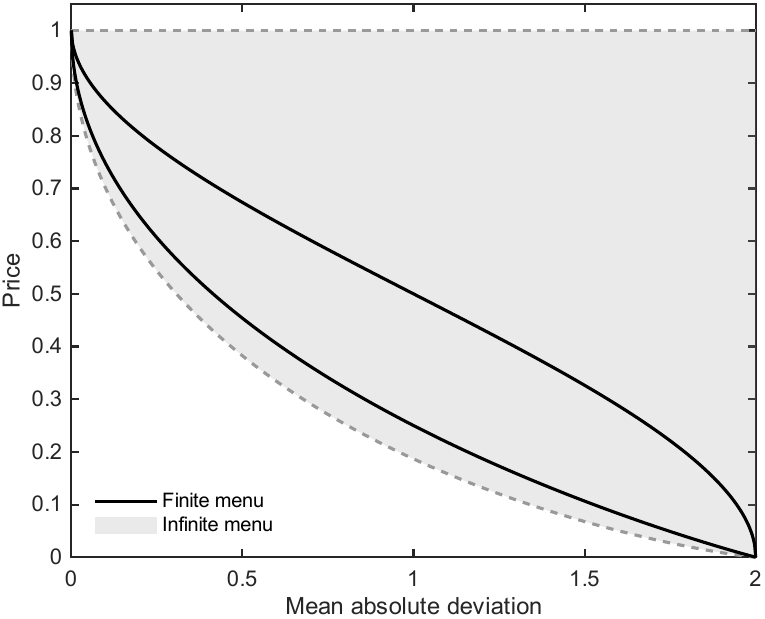}
  \end{subfigure}
 
  \vspace{1ex}
 
  \begin{subfigure}[b]{0.48\linewidth}
    \centering
    \caption*{$N = 4$}
    \includegraphics[width=\linewidth]{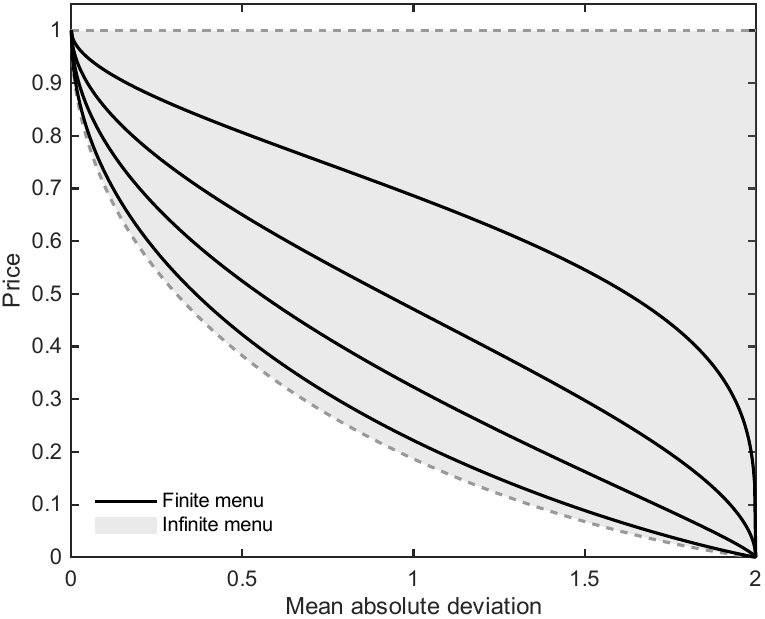}
  \end{subfigure}\hfill
  \begin{subfigure}[b]{0.48\linewidth}
    \centering
    \caption*{$N = 8$}
    \includegraphics[width=\linewidth]{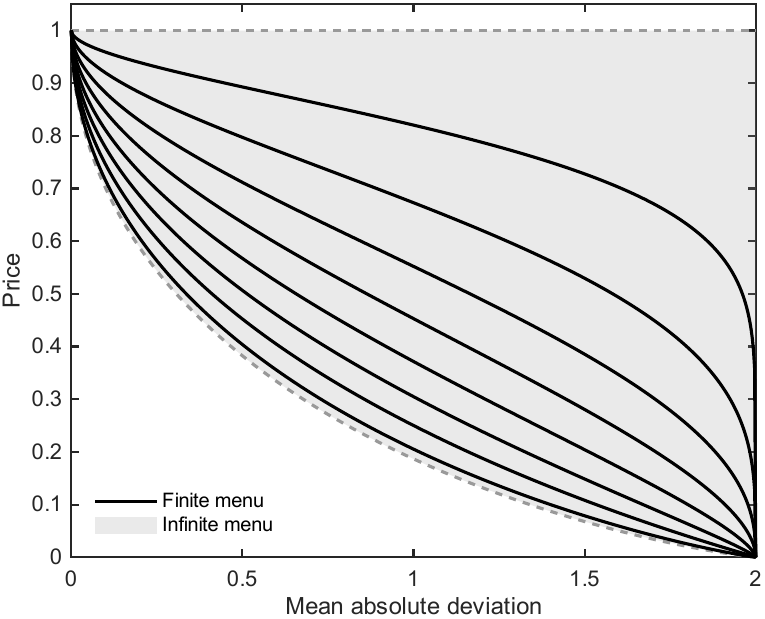}
  \end{subfigure}
  \caption{Optimal finite-menu prices for $\cP = \cP_d(\mu,d)$ and different values of $N$. The finite-menu prices are compared to the infinite-menu prices, as functions of the MAD with $\mu=1$.}
  \label{fig:prices_mad}
\end{figure}

Figure~\ref{fig:prices_mad} illustrates Proposition~\ref{prop:prices_decrease_MAD}. For every menu size, each price of the optimal mean--MAD menu decreases as the MAD increases, in contrast with the mean--variance menus in Figure~\ref{fig:prices_var}, whose higher prices turn upward as dispersion grows. All prices remain below the mean, and as the menu grows, they fill the support of the fully randomized mechanism.

\end{appendix}
\end{document}